\documentclass[preprint,12pt,authoryear] {elsarticle}

\usepackage{epsfig,graphics,graphicx,amsmath,amssymb,amsfonts}
\usepackage{forest}
\usepackage{arydshln}
\usepackage{stmaryrd}
\usepackage{subfigure}
\usepackage[colorlinks,bookmarksopen,bookmarksnumbered,citecolor=red,urlcolor=red]{hyperref}
\usepackage[percent]{overpic}
\usepackage[makeroom]{cancel}
\usepackage{epstopdf}
\usepackage{hyperref}
\usepackage{natbib}
\setcitestyle{numbers,square}
\def\d{\partial}

\usepackage[margin=1in]{geometry}
\usepackage{color}
\usepackage[title]{appendix}

\renewcommand{\tilde}{\widetilde}
\newcommand{\eps}{\varepsilon}
\newtheorem{remark}{Remark}
\numberwithin{remark}{section}

\newcommand{\bm}[1]{\mbox{\boldmath{$#1$}}}

\def\0{{\bm 0}}

\journal{arXiv}
\begin{document}

\begin{frontmatter}



\title{Extended one-dimensional reduced model for blood flow within a stenotic artery} 


\author[l1]{Suncica Canic} 
\affiliation[l1]{organization={Department of Mathematics, University of California Bekeley},
            city={Berkeley},
            postcode={94720}, 
            state={CA},
            country={USA}}

\author[l2]{Shihan Guo} 
\affiliation[l2]{organization={School of Mathematical Sciences, East China Normal University},
            city={Shanghai},
            postcode={200241}, 
            state={Shanghai},
            country={China}}

\author[l3]{Yifan Wang} 
\affiliation[l3]{organization={Department of Mathematics and Statistics, Texas Tech University},
            city={Lubbock},
            postcode={79409}, 
            state={TX},
            country={USA}}

\author[l2]{Xiaohe Yue} 
            
\author[l4]{Haibiao Zheng} 
\affiliation[l4]{organization={School of Mathematical Sciences, Ministry of Education Key Laboratory of Mathematics and Engineering Applications, Shanghai Key Laboratory of PMMP,  East China Normal University},
            city={Shanghai},
            postcode={200241}, 
            state={Shanghai},
            country={China}}

\begin{abstract}
In this paper, we introduce an adapted one-dimensional ($1D$) reduced model aimed at analyzing blood flow within stenosed arteries. Differing from the prevailing $1D$ model \cite{Formaggia2003, Sherwin2003_2, Sherwin2003, Quarteroni2004, 10.1007/978-3-642-56288-4_10}, our approach incorporates the variable radius of the blood vessel. Our methodology begins with the non-dimensionalization of the Navier-Stokes equations for axially symmetric flow in cylindrical coordinates and then derives the extended $1D$ reduced model, by making additional adjustments to accommodate the effects of variable radii of the vessel along the longitudinal direction. Additionally, we propose a method to extract radial velocity information from the $1D$ results during post-processing, enabling the generation of two-dimensional ($2D$) velocity data. We validate our model by conducting numerical simulations of blood flow through stenotic arteries with varying severities, ranging from 23\% to 50\%. The results were compared to those from the established $1D$ model and a full three-dimensional ($3D$) simulation, highlighting the potential and importance of this model for arteries with variable radius. All the code used to generate the results presented in the paper is available at \url{https://github.com/qcutexu/Extended-1D-AQ-system.git}.
\end{abstract}

\begin{graphicalabstract}
The extended $1D$ model in conservative form is expressed as:
\begin{equation}
\left\{
\begin{aligned}
&\frac{\partial A}{\partial t}
+ \frac{\partial Q}{\partial z} = 0,\nonumber\\
&\frac{\partial Q}{\partial t}
+ \frac{\partial}{\partial z}\left[(\alpha{+\alpha_c}) \frac{Q^{2}}{A }\right]
+\frac{A}{\rho_{f}} \frac{\partial p}{\partial z}
=-2(\gamma+2) \nu \frac{Q}{A}
+\frac{Q^2}{A}\frac{\partial\alpha_c}{\partial z},\nonumber\\
&p=
p_{\text {ext }}
+ \frac{h E}{R_{0}^2\left(1-\sigma^{2}\right)}\left({R}-{R_{0}}\right)
+ (\gamma+2) \rho_{f} \nu \frac{Q}{ A } \frac{\partial \ln R_{0}}{\partial z}.\nonumber
\end{aligned}
\right.
\end{equation}
here, $\gamma=-\frac{\alpha-2}{\alpha-1}$, and $\alpha$ is defined by: 
\begin{equation}
    \alpha_c = -\frac{2}{35}\left(\frac{\partial R_0}{\partial z}\right)^2.
    \nonumber
\end{equation}
When compared to the existing $1D$ model:
\begin{equation}
\left\{
\begin{aligned}
&\frac{\partial A}{\partial t}
+ \frac{\partial Q}{\partial z} = 0,\nonumber\\
&\frac{\partial Q}{\partial t}
+ \frac{\partial}{\partial z}\left(\alpha \frac{Q^{2}}{A }\right)
+\frac{A}{\rho_{f}} \frac{\partial p}{\partial z}
=-2(\gamma+2) \nu \frac{Q}{A},\nonumber\\
&p=
p_{\text {ext }}
+ \frac{h E}{R_{0}^2\left(1-\sigma^{2}\right)}\left({R}-{R_{0}}\right).\nonumber
\end{aligned}
\right.
\end{equation}
The extended model incorporates a few additional terms which account for spatial variations in the arterial radius, making it more suitable for modeling cases involving stenotic arteries. This added complexity allows the model to accurately handle flow behavior and pressure variations in the presence of arterial narrowing (stenosis) as shown in Figure \ref{fig_1}, where changes in the geometry can significantly impact blood flow dynamics. 
\vspace{0.3in}
\begin{figure}[ht]
  \centering
 \begin{overpic}[width=0.88\textwidth,grid=false,tics=10]{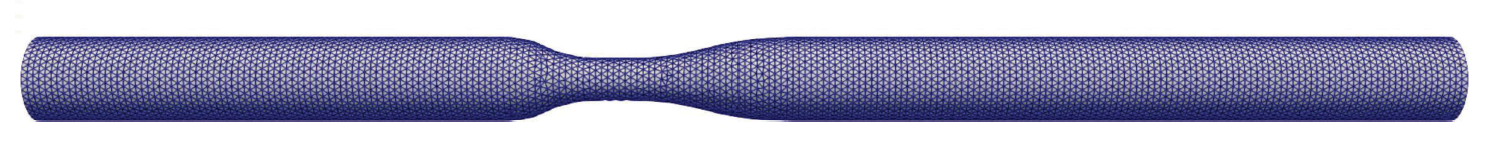}
      \put(-5,4){{Inlet}}
      \put(99,4){{Outlet}}
    \put(10,-1){{50\% stenosis:   ~ $R_0 =R_{max}-0.5R_{max}*e^{-50( z-3.4+0.95e^{-0.5*(z-2.5)^2 })^4}$}}
    \end{overpic}
  \caption{
The $3D$ mesh representing the domain of a stenotic artery with a 50\% narrowing.}
\label{fig_1}
\end{figure}

\end{graphicalabstract}

\begin{highlights}
\item Introduced an extended $1D$ reduced model for blood flow in stenotic arteries, incorporating the variable radius of the blood vessel, unlike the established $1D$ models.
\item The extended $1D$ reduced model successfully captures the accurate solution, matching that of the full $3D$ model, whereas the established $1D$ models fail to produce the correct results. 
\item This work also developed a method to extract radial velocity from $1D$ results, providing $2D$ velocity data for improved analysis.
\end{highlights}


\end{frontmatter}

\section{Introduction}
The complex nature of the cardiovascular system demands accurate and efficient models to simulate blood flow dynamics for research and clinical applications. While three-dimensional ($3D$) models provide detailed insights into the vascular system \cite{Janela2010}, they are often computationally intensive and challenging to implement in large-scale simulations or provide real-time feedback, making it difficult to provide real-time feedback for clinical decision making. Therefore, there is a growing interest in developing reduced-order models, such as one-dimensional ($1D$) and two-dimensional ($2D$) blood flow models, as referenced in \cite{Formaggia2003,Sherwin2003_2, Quarteroni2004, 10.1007/978-3-642-56288-4_10, Figueroa2005, Boujena2014, Arthur2017, Quarteroni2016, Saito2011, Tambaca2005}. These models offer a practical compromise between computational efficiency and physiological fidelity.

The fundamental principles governing the $1D$ blood flow model are based on the conservation laws of mass and momentum. By simplifying the cardiovascular system and representing blood vessels as cylindrical tubes, this model averages the flow variables across the vessel's cross-section. This simplification leads to a system of hyperbolic partial differential equations that describes the spatio-temporal evolution of key variables such as average flow rate and cross-sectional area. Once the system is solved, blood flow velocity and pressure can be derived from key variables. This dimensional reduction enables the simulation of large vascular networks within a feasible computational time while capturing essential hemodynamics, including pressure and flow wave propagation, reflections, and interactions between different vascular segments.
Despite its simplifications, the $1D$ model has been proven to accurately predict the physiological responses of the cardiovascular system under various conditions. It has been used to investigate the impact of arterial stiffness and stenosis on blood pressure wave \cite{Sherwin2003_2, Jin2021, Duanmu2019}. Its capacity to integrate with larger-scale models and experimental data makes it highly applicable in personalized medicine, facilitating patient-specific simulations that can inform clinical decisions and treatment planning. \cite{Saito2011, Mynard2008, YIN201966, Epstein2015, Padmos2019, Caforio2022}.

\begin{figure}[htbp]
	\center
\includegraphics[width=0.8\columnwidth]{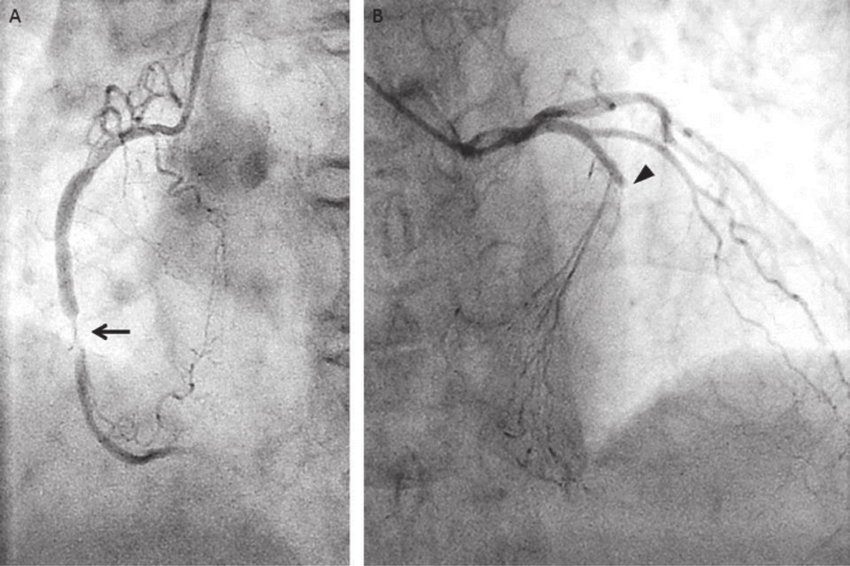}
\caption{Angiography of the stenosis with a steep change \cite{Shirai2020}. (A): Severe stenosis is seen in the mid-portion of the right coronary artery. (B) Total occlusion is observed in the left anterior descending coronary artery.}
\label{angiography}
\end{figure}

Although the existing $1D$ model offers computational efficiency and simplicity, its limitations, such as the assumptions of cylindrical vessel geometry and uniform flow distribution, hinder its applicability in modeling stenosed arteries and other complex vascular conditions. To address these limitations, we propose an extended $1D$ reduced model for studying blood flow in an artery with a variable radius. This modification is motivated by the necessity in the case of analyzing blood flow through stenosis regions formed by atherosclerotic plaque or other substances (vascular plasticity) within the walls of arteries, where the artery's radius varies along the longitudinal direction, as illustrated in Figure \ref{angiography}. Unlike the well-established $1D$ model introduced in \cite{Formaggia2003,  Sherwin2003, Figueroa2005}, we make adjustments to the model that depend on both the artery's radius and its derivative. Furthermore, we propose a method to extract radial velocity information from our $1D$ results during the post-processing phase, allowing us to derive $2D$ velocity data from the $1D$ model.

This manuscript follows the following structure:
Section 2 provides a detailed derivation of the extended $1D$ model, considering variable arterial radius.
In Section 3, we outline our numerical approach for solving this system, based on the discontinuous Galerkin (DG) method.
Section 4 presents the results of our numerical simulations, highlighting how the proposed model correctly captures the flow dynamics in arteries with stenosis, where existing $1D$ models fail to yield accurate outcomes.
Section 5 provides a summary of the key contributions of this study, emphasizes the importance of the proposed model, and suggests potential applications along with future research directions.

\section{Derivation of the extended $1D$ reduced model for blood flow with variable radius}
In our model, we assume blood behaves as an incompressible, Newtonian viscous fluid, and the flow inside the artery is axially symmetric, allowing us to disregard angular effects. We consider the arterial wall to be homogeneous and isotropic, with the thickness of the wall $h$ and its radius $R_0$ significantly smaller than the length of the arterial segment $L$. We model the mechanical properties of arterial walls using the thin membrane structure, specifically the Koiter shell model. The subsequent sections will report both the original $3D$ equations and the derivation of the reduced model in details.

\subsection{The Navier-Stokes equation in cylindrical coordinates assuming variable radius}
We initiate our derivation by examining the motion of an incompressible Newtonian fluid within an axial symmetric cylindrical domain with a variable radius. The length of the domain is denoted by $L$. For a given smooth function $R_0:[0, L] \rightarrow \mathbb{R}$, the radius of the cylinder at $z \in[0, L]$ is denoted by $R_0(z)$. The reference domain is now defined by:
$$
\Omega=\left\{x=(r \cos \theta, r \sin \theta, z) \in \mathbb{R}^3: r \in(0, R(z)), \theta \in(0,2 \pi), z \in(0, L)\right\},
$$
and its lateral boundary is given by
$$
\Gamma=\left\{x=(R_0(z) \cos \theta, R_0(z) \sin \theta, z) \in \mathbb{R}^3: \theta \in(0,2 \pi), z \in(0, L)\right\}.
$$
The Navier-Stokes (NS) equations governing axially symmetric flow in $\Omega$ are expressed as:
\begin{equation}
\left\{~
\begin{aligned}
&\frac{\partial u_{r}}{\partial r}+\frac{\partial u_{z}}{\partial z}+\frac{u_{r}}{r} =0,\\
&\rho_{f}\left(\frac{\partial u_{r}}{\partial t}+u_{r} \frac{\partial u_{r}}{\partial r}+u_{z} \frac{\partial u_{r}}{\partial z}\right)-\mu_{f}\left(\frac{\partial^{2} u_{r}}{\partial r^{2}}+\frac{\partial^{2} u_{r}}{\partial z^{2}}+\frac{1}{r} \frac{\partial u_{r}}{\partial r}-\frac{u_{r}}{r^{2}}\right)+\frac{\partial p}{\partial r} =0, \\
&\rho_{f}\left(\frac{\partial u_{z}}{\partial t}+u_{r} \frac{\partial u_{z}}{\partial r}+u_{z} \frac{\partial u_{z}}{\partial z}\right)-\mu_{f}\left(\frac{\partial^{2} u_{z}}{\partial r^{2}}+\frac{\partial^{2} u_{z}}{\partial z^{2}}+\frac{1}{r} \frac{\partial u_{z}}{\partial r}\right)+\frac{\partial p}{\partial z} =0. 
\end{aligned}
\right.
\label{NS}
\end{equation}
Here, $u_r$, $u_z$, and $p$ represent the radial velocity, axial velocity, and pressure respectively. $\mu_f$ and $\rho_f$ denote the dynamic viscosity and density of the blood respectively. It is worth noting that Equ. (\ref{NS}.a), concerning mass conservation, is not stated in conservative form due to the assumption of longitudinal variation in radius caused by stenosis.

\subsection{Derivation of the system in non-dimensional form}\label{sec.derivatives}
To derive the non-dimensional form of the NS system (\ref{NS}), we introduce the following non-dimensional variables:
$$
r=R_{0}(\tilde{z}) \tilde{r},~ z = L \tilde{z},~ u_{z}=U_{z} \tilde{u}_{z},~ u_{r}=U_{r}(\tilde{z}) \tilde{u}_{r},~ t=T \tilde{t}=\frac{L}{U_{z}} \tilde{t},~ p=\rho_{f} U_{z}^{2} \tilde{p},
$$
where $R_0(\tilde{z})$ represents the characteristic radius, $L$ denotes the characteristic length, $U_r(\tilde{z})$ is the characteristic radial velocity, $U_z$ stands for the characteristic axial velocity, and $T$ indicates the characteristic time.

To account for the stenosis (narrowing) in the artery, we consider $R_0(\tilde{z})$ to vary with the axial coordinate $\tilde{z}$. For consistency in the scaling, we require that the ratio $\frac{U_z ~R_0(\tilde{z})}{L ~U_r} = 1$, which leads to the assumption that $U_r = \frac{U_z R_0(\tilde{z})}{L}$, meaning that the radial velocity also varies with $\tilde{z}$. Additionally, we assume the radial and axial length scale satisfies the relation:
$$
\frac{R_{0}(\tilde{z})}{L}=\varepsilon(\tilde{z}) \ll 1.
$$
For the time scale, we also have $t=\frac{L}{U_{z}} \tilde{t}=\frac{R_{0}(\tilde{z})}{U_{r}(\tilde{z})} \tilde{t}$. It is important to note that the coordinate mapping between $\tilde{r}, \tilde{z}$ and $r, z$ implies that $r$ and $z$ are not independent variables. For the unknown variables $u_r$, $u_z$, and $p$, the following relations hold.

\hspace{0.25in}
$
\begin{forest}
  for tree={
    math content,
    calign=fixed edge angles,
  },
  before typesetting nodes={
    for descendants={
      delay={content/.wrap value={\strut #1}},
      inner sep=1.5pt,
      where={
        >OOw+n< {n}{!u.n children}{(#1+1)/2}
      }{
        edge label/.process={
          OOw2 {content}{!u.content} {node [midway, left, anchor=base east]{ }
          }
        }
      }{
        edge label/.process={
          OOw2 {content}{!u.content} {node [midway, right, anchor=base west] { }
          }
        }
      },
    },
  },
  [u_r
   [\tilde{u}_r
      [\tilde{r}
        [r]
        [z]
      ]
      [\tilde{z}
        [z]
      ]
      [\tilde{t}
        [t]
      ]
    ]
    [U_r(\tilde{z})
      [\tilde{z}
        [z]
      ]
    ]
  ]  
\end{forest}
$
\hspace{0.5in}
$
\begin{forest}
  for tree={
    math content,
    calign=fixed edge angles,
  },
  before typesetting nodes={
    for descendants={
      delay={content/.wrap value={\strut #1}},
      inner sep=1.5pt,
      where={
        >OOw+n< {n}{!u.n children}{(#1+1)/2}
      }{
        edge label/.process={
          OOw2 {content}{!u.content} {node [midway, left, anchor=base east] { } 
          }
        }
      }{
        edge label/.process={
          OOw2 {content}{!u.content} {node [midway, right, anchor=base west] { } 
          }
        }
      },
    },
  },
  [u_z
   [\tilde{u}_z
      [\tilde{r}
        [r]
        [z]
      ]
      [\tilde{z}
        [z]
      ]
      [\tilde{t}
        [t]
      ]
    ]
    [U_z
    ]
  ]  
\end{forest}
$
\hspace{0.5in}
$
\begin{forest}
  for tree={
    math content,
    calign=fixed edge angles,
  },
  before typesetting nodes={
    for descendants={
      delay={content/.wrap value={\strut #1}},
      inner sep=1.5pt,
      where={
        >OOw+n< {n}{!u.n children}{(#1+1)/2}
      }{
        edge label/.process={
          OOw2 {content}{!u.content} {node [midway, left, anchor=base east]{ } 
          }
        }
      }{
        edge label/.process={
          OOw2 {content}{!u.content} {node [midway, right, anchor=base west]{ } 
          }
        }
      },
    },
  },
   [p
   [\tilde{p}
      [\tilde{r}
        [r]
        [z]
      ]
      [\tilde{z}
        [z]
      ]
      [\tilde{t}
        [t]
      ]
    ]
    [U_z]
    ]
\end{forest}
$
\vspace{0.1in}
Utilizing the chain rule, we derive the following expressions for the first-order derivative terms appearing in Eqs.~(\ref{NS}):
$$
\left\{
\begin{aligned}
\frac{\partial u_r}{\partial r} &= \frac{\partial u_r}{\partial \tilde{u}_r} \frac{\partial \tilde{u}_r}{\partial \tilde{r}}\frac{\partial \tilde{r}}{\partial r} = \frac{U_r(\tilde{z})}{R_0(\tilde{z})}\frac{\partial \tilde{u}_r}{\partial \tilde{r}}, \\
\frac{\partial u_r}{\partial z} &= \frac{\partial u_r}{\partial \tilde{u}_r}\left( \frac{\partial \tilde{u}_r}{\partial \tilde{r}}\frac{\partial \tilde{r}}{\partial z} + \frac{\partial \tilde{u}_r}{\partial \tilde{z}}\frac{\partial \tilde{z}}{\partial z} \right)+\frac{\partial u_r}{\partial U_r(\tilde{z})}\frac{\partial U_r(\tilde{z})}{\partial\tilde{z}}\frac{\partial\tilde{z}}{\partial z} \\
&= -\frac{U_r(\tilde{z})\tilde{r}}{L}\frac{\partial lnR_0(\tilde{z})}{\partial\tilde{z}}\frac{\partial\tilde{u}_r}{\partial\tilde{r}}+\frac{U_r(\tilde{z})}{L}\frac{\partial\tilde{u}_r}{\partial\tilde{z}} + \frac{\tilde{u}_r}{L}\frac{\partial U_r(\tilde{z})}{\partial\tilde{z}},   \\ 
\frac{\partial u_r}{\partial t} &=\frac{\partial u_r}{\partial\tilde{u}_r}\frac{\partial\tilde{u}_r}{\partial\tilde{t}}\frac{\partial\tilde{t}}{\partial t} = \frac{U_r(\tilde{z})U_z}{L}\frac{\partial\tilde{u}_r}{\partial\tilde{t}}.
\end{aligned}
\right.
$$

$$
\left\{
\begin{aligned}
\frac{\partial u_z}{\partial r} &= \frac{\partial u_z}{\partial\tilde{u}_z}\frac{\partial\tilde{u}_z}{\partial\tilde{r}}\frac{\partial\tilde{r}}{\partial r} = \frac{U_z}{R_0(\tilde{z})}\frac{\partial\tilde{u}_z}{\partial\tilde{r}},    \\             
\frac{\partial u_z}{\partial z} &= \frac{\partial u_z}{\partial \tilde{u}_z}\left( \frac{\partial\tilde{u}_z}{\partial \tilde{z}}\frac{\partial \tilde{z}}{\partial z} + \frac{\partial \tilde{u}_z}{\partial \tilde{r}}\frac{\partial \tilde{r}}{\partial z} \right) 
= \frac{U_z}{L}\frac{\partial \tilde{u}_z}{\partial \tilde{z}} - \frac{U_z\tilde{r}}{L}\frac{\partial ln(R_0(\tilde{z}))}{\partial \tilde{z}}\frac{\partial \tilde{u}_z}{\partial \tilde{r}},\\
\frac{\partial u_z}{\partial t} &= \frac{\partial u_z}{\partial\tilde{u}_z}\frac{\partial\tilde{u}_z}{\partial\tilde{t}}\frac{\partial\tilde{t}}{\partial t} = \frac{U_z^2}{L}\frac{\partial\tilde{u}_z}{\partial\tilde{t}}.
\end{aligned}
\right.
$$

$$
\left\{
\begin{aligned}
\frac{\partial p}{\partial r} &= \frac{\partial p}{\partial\tilde{p}}\frac{\partial\tilde{p}}{\partial\tilde{r}}\frac{\partial\tilde{r}}{\partial r} = \frac{\rho_fU_z^2}{R_0(\tilde{z})}\frac{\partial\tilde{p}}{\partial\tilde{r}},  \\
\frac{\partial p}{\partial z} &= \frac{\partial p}{\partial\tilde{p}}\left( \frac{\partial\tilde{p}}{\partial\tilde{z}}\frac{\partial\tilde{z}}{\partial z} + \frac{\partial\tilde{p}}{\partial\tilde{r}}\frac{\partial\tilde{r}}{\partial z} \right) = \frac{\rho_fU_z^2}{L}\frac{\partial\tilde{p}}{\partial\tilde{z}} - \frac{\rho_fU_z^2\tilde{r}}{L}\frac{\partial ln(R_0(\tilde{z}))}{\partial\tilde{z}}\frac{\partial\tilde{p}}{\partial\tilde{r}}.    
\end{aligned}
\right.
$$
Regarding the second-order derivative terms present in Eqs.~(\ref{NS}), we have:
$$
\left\{
\begin{aligned}
\frac{\partial^2 u_r}{\partial r^2} =&~ \frac{\partial}{\partial r}\left( \frac{U_r(\tilde{z})}{R_0(\tilde{z})}\frac{\partial \tilde{u}_r}{\partial \tilde{r}} \right) = \frac{U_r(\tilde{z})}{R_0^2(\tilde{z})}\frac{\partial^2\tilde{u}_r}{\partial\tilde{r}^2},\\
\frac{\partial^2 u_r}{\partial z^2} =&~ \frac{\partial}{\partial z}\left( -\frac{U_r(\tilde{z})\tilde{r}}{L}\frac{\partial ln R_0(\tilde{z})}{\partial\tilde{z}}\frac{\partial\tilde{u}_r}{\partial\tilde{r}}+\frac{U_r(\tilde{z})}{L}\frac{\partial\tilde{u}_r}{\partial\tilde{z}} + \frac{\tilde{u}_r}{L}\frac{\partial U_r(\tilde{z})}{\partial\tilde{z}}  \right) \\
=& -\frac{\tilde{r}}{L^2}\frac{\partial U_r(\tilde{z})}{\partial\tilde{z}}\frac{\partial ln R_0(\tilde{z})}{\partial\tilde{z}}\frac{\partial\tilde{u}_r}{\partial\tilde{r}} 
+ \frac{U_r(\tilde{z})\tilde{r}}{L^2}\frac{\partial ln R_0(\tilde{z})}{\partial\tilde{z}}\frac{\partial ln R_0(\tilde{z})}{\partial\tilde{z}}\frac{\partial\tilde{u}_r}{\partial\tilde{r}} 
- \frac{U_r(\tilde{z})\tilde{r}}{L^2}\frac{\partial^2 ln R_0(\tilde{z})}{\partial\tilde{z}^2}\frac{\partial\tilde{u}_r}{\partial\tilde{r}}\\
& - \frac{U_r(\tilde{z})\tilde{r}}{L^2}\frac{\partial ln R_0(\tilde{z})}{\partial\tilde{z}}\frac{\partial^2\tilde{u}_r}{\partial\tilde{r}\partial\tilde{z}} + \frac{U_r(\tilde{z})\tilde{r}^2}{L^2}\frac{\partial ln R_0(\tilde{z})}{\partial\tilde{z}}\frac{\partial ln R_0(\tilde{z})}{\partial\tilde{z}}\frac{\partial^2\tilde{u}_r}{\partial\tilde{r}^2} \\
&+\frac{1}{L^2}\frac{\partial U_r(\tilde{z})}{\partial\tilde{z}}\frac{\partial\tilde{u}_r}{\partial\tilde{z}}+\frac{U_r(\tilde{z})}{L^2}\frac{\partial^2\tilde{u}_r}{\partial\tilde{z}^2} - \frac{U_r(\tilde{z})\tilde{r}}{L^2}\frac{\partial ln R_0(\tilde{z})}{\partial\tilde{z}}\frac{\partial^2\tilde{u}_r}{\partial\tilde{z}\partial\tilde{r}} \\
&+ \frac{\tilde{u}_r}{L^2}\frac{\partial^2 U_r(\tilde{z})}{\partial\tilde{z}^2} + \frac{1}{L^2}\frac{\partial U_r(\tilde{z})}{\partial\tilde{z}}\frac{\partial\tilde{u}_r}{\partial\tilde{z}} - \frac{\tilde{r}}{L^2}\frac{\partial U_r(\tilde{z})}{\partial\tilde{z}}\frac{\partial ln R_0(\tilde{z})}{\partial\tilde{z}}\frac{\partial\tilde{u}_r}{\partial\tilde{r}},\\
\frac{\partial^2 u_z}{\partial r^2} 
= &~ \frac{\partial}{\partial r}\left( \frac{U_z}{R_0(\tilde{z})}\frac{\partial\tilde{u}_z}{\partial\tilde{r}} \right) 
= \frac{U_z}{R_0^2(\tilde{z})}\frac{\partial^2\tilde{u}_z}{\partial\tilde{r}^2},\\
\frac{\partial^2 u_z}{\partial z^2} =&~ \frac{\partial}{\partial z}\left(\frac{U_z}{L}\frac{\partial \tilde{u}_z}{\partial \tilde{z}} - \frac{U_z\tilde{r}}{L}\frac{\partial ln R_0(\tilde{z})}{\partial \tilde{z}}\frac{\partial \tilde{u}_z}{\partial \tilde{r}} \right) \\
=&~ \frac{U_z}{L^2}\frac{\partial^2\tilde{u}_z}{\partial\tilde{z}^2} 
- \frac{U_z\tilde{r}}{L^2}\frac{\partial ln R_0(\tilde{z})}{\partial\tilde{z}}\frac{\partial^2\tilde{u}_z}{\partial\tilde{z}\partial\tilde{r}} 
+ \frac{U_z\tilde{r}}{L^2}\frac{\partial ln R_0(\tilde{z})}{\partial\tilde{z}}\frac{\partial ln R_0(\tilde{z})}{\partial\tilde{z}}\frac{\partial\tilde{u}_z}{\partial\tilde{r}} - \frac{U_z\tilde{r}}{L^2}\frac{\partial^2ln R_0(\tilde{z})}{\partial\tilde{z}^2}\frac{\partial\tilde{u}_z}{\partial\tilde{r}}\\
&+ \frac{U_z\tilde{r}^2}{L^2}\frac{\partial ln R_0(\tilde{z})}{\partial\tilde{z}}\frac{\partial ln R_0(\tilde{z})}{\partial\tilde{z}}\frac{\partial^2\tilde{u}_z}{\partial\tilde{r}^2} - \frac{U_z\tilde{r}}{L^2}\frac{\partial^2\tilde{u}_z}{\partial\tilde{r}\partial\tilde{z}},   \\
\frac{u_r}{r} =&~ \frac{U_r}{R_0(\tilde{z})}\frac{\tilde{u}_r}{\tilde{r}}.
\end{aligned}
\right.
$$


\subsection{The conservation of mass in nondimensional form}
Substituting the derived expressions in Sec.~\ref{sec.derivatives} into the Equ.~(\ref{NS}.a), we have:
\begin{equation}
\frac{\partial \tilde{u}_r}{\partial \tilde{r}} + \frac{U_zR_0(\tilde{z})}{U_r(\tilde{z})L}\frac{\partial \tilde{u}_z}{\partial \tilde{z}} + \frac{\tilde{u}_r}{\tilde{r}} - \frac{U_z\tilde{r}}{U_r(\tilde{z})L}\frac{\partial R_0(\tilde{z})}{\partial \tilde{z}}\frac{\tilde{u}_z}{\partial \tilde{r}} = 0
\end{equation}
Multiplying the aforementioned equation by $\tilde{r}$ and observing that $\frac{1}{R_0(\tilde z)}\frac{\partial R_0(\tilde z)}{\partial \tilde{z}}=\frac{\partial \ln R_0(\tilde z)}{\partial \tilde{z}}$, along with acknowledging the condition $\frac{U_{z} R_0(\tilde{z})}{U_{r}(\tilde{z}) L}=1$, we derive the following equation:
	\begin{equation}
		\frac{\partial (\tilde{r} \tilde{u}_{r})}{\partial \tilde{r}} + \frac{\partial (\tilde{r} \tilde{u}_{z})}{\partial \tilde{z}} 
 		- \tilde{r}^2\frac{\partial \ln R_{0}(\tilde z)}{\partial \tilde{z}}\frac{\partial \tilde{u}_z}{\partial \tilde{r}}
 		=0, 
		\label{mass_conser}
	\end{equation}
Noting that the emergence of the third term is a result of the variation of $R_0(\tilde{z})$ with respect to $\tilde{z}$.

\subsection{The radial momentum equation in nondimensional form}
Substitute in the terms we derived in Sec.~\ref{sec.derivatives} for Equ.~(\ref{NS}.b), we have:
$$
\begin{aligned}
\frac{U_r(\tilde{z})U_z}{L}\frac{\partial\tilde{u}_r}{\partial\tilde{t}}
+ U_r(\tilde{z})\tilde{u}_r\frac{U_r(\tilde{z})}{R_0(\tilde{z})}\frac{\partial \tilde{u}_r}{\partial \tilde{r}}
+ U_z\tilde{u}_z\left(-\frac{U_r(\tilde{z})\tilde{r}}{L}\frac{\partial ln R_0(\tilde{z})}{\partial\tilde{z}}\frac{\partial\tilde{u}_r}{\partial\tilde{r}}
+\frac{U_r(\tilde{z})}{L}\frac{\partial\tilde{u}_r}{\partial\tilde{z}} + \frac{\tilde{u}_r}{L}\frac{\partial U_r(\tilde{z})}{\partial\tilde{z}}\right)\\
+\frac{U_z^2}{R_0(\tilde{z})}\frac{\partial\tilde{p}}{\partial\tilde{r}}  
= \frac{\mu_f}{\rho_f}\Bigg\{\frac{U_r(\tilde{z})}{R_0^2(\tilde{z})}\frac{\partial^2\tilde{u}_r}{\partial\tilde{r}^2}
-\frac{\tilde{r}}{L^2}\frac{\partial U_r(\tilde{z})}{\partial\tilde{z}}\frac{\partial ln R_0(\tilde{z})}{\partial\tilde{z}}\frac{\partial\tilde{u}_r}{\partial\tilde{r}} 
+ \frac{U_r(\tilde{z})\tilde{r}}{L^2}\frac{\partial ln R_0(\tilde{z})}{\partial\tilde{z}}\frac{\partial ln R_0(\tilde{z})}{\partial\tilde{z}}\frac{\partial\tilde{u}_r}{\partial\tilde{r}}\\ 
- \frac{U_r(\tilde{z})\tilde{r}}{L^2}\frac{\partial^2 ln R_0(\tilde{z})}{\partial\tilde{z}^2}\frac{\partial\tilde{u}_r}{\partial\tilde{r}} 
- \frac{U_r(\tilde{z})\tilde{r}}{L^2}\frac{\partial ln R_0(\tilde{z})}{\partial\tilde{z}}\frac{\partial^2\tilde{u}_r}{\partial\tilde{r}\partial\tilde{z}} 
+ \frac{U_r(\tilde{z})\tilde{r}^2}{L^2}\frac{\partial ln R_0(\tilde{z})}{\partial\tilde{z}}\frac{\partial ln R_0(\tilde{z})}{\partial\tilde{z}}\frac{\partial^2\tilde{u}_r}{\partial\tilde{r}^2} \\
+\frac{1}{L^2}\frac{\partial U_r(\tilde{z})}{\partial\tilde{z}}\frac{\partial\tilde{u}_r}{\partial\tilde{z}}+\frac{U_r(\tilde{z})}{L^2}\frac{\partial^2\tilde{u}_r}{\partial\tilde{z}^2} 
- \frac{U_r(\tilde{z})\tilde{r}}{L^2}\frac{\partial ln R_0(\tilde{z})}{\partial\tilde{z}}\frac{\partial^2\tilde{u}_r}{\partial\tilde{z}\partial\tilde{r}} \\
+ \frac{\tilde{u}_r}{L^2}\frac{\partial^2 U_r(\tilde{z})}{\partial\tilde{z}^2} 
+ \frac{1}{L^2}\frac{\partial U_r(\tilde{z})}{\partial\tilde{z}}\frac{\partial\tilde{u}_r}{\partial\tilde{z}} 
- \frac{\tilde{r}}{L^2}\frac{\partial U_r(\tilde{z})}{\partial\tilde{z}}\frac{\partial ln R_0(\tilde{z})}{\partial\tilde{z}}\frac{\partial \tilde{u}_r}{\partial\tilde{r}}\\
+ \frac{1}{R_0(\tilde{z})\tilde{r}}\frac{U_r(\tilde{z})}{R_0(\tilde{z})}\frac{\partial \tilde{u}_r}{\partial \tilde{r}} 
- \frac{U_r(\tilde{z})\tilde{u}_r}{R_0(\tilde{z})^2\tilde{r}^2}
\Bigg\}.
\end{aligned}
$$
We first divide the equation above by $U_{z}^{2}$, and then multiply it by $R_{0}(\tilde{z})$, resulting in:
$$
\begin{aligned}
\frac{U_r(\tilde{z})R_0(\tilde{z})}{LU_z}\frac{\partial\tilde{u}_r}{\partial\tilde{t}}
+ \tilde{u}_r\frac{U_r(\tilde{z})^2}{U_z^2}\frac{\partial \tilde{u}_r}{\partial \tilde{r}}
- \tilde{u}_z\tilde{r}\frac{U_r(\tilde{z})R_0(\tilde{z})}{LU_z}\frac{\partial ln R_0(\tilde{z})}{\partial\tilde{z}}\frac{\partial\tilde{u}_r}{\partial\tilde{r}}
+ \tilde{u}_z\frac{U_r(\tilde{z})R_0(\tilde{z})}{LU_z}\frac{\partial\tilde{u}_r}{\partial\tilde{z}}\\
+ \frac{R_0(\tilde{z})\tilde{u}_r\tilde{u}_z}{LU_z}\frac{\partial U_r(\tilde{z})}{\partial\tilde{z}}
+ \frac{\partial\tilde{p}}{\partial\tilde{r}} 
=\frac{\mu_f L}{\rho_f R_0(\tilde{z})^2U_z}\Bigg\{
	\frac{R_0(\tilde{z})U_r(\tilde{z})}{LU_z}\frac{\partial^2\tilde{u}_r}{\partial\tilde{r}^2}
-\frac{R_0(\tilde{z})^3\tilde{r}}{L^3U_z}\frac{\partial U_r(\tilde{z})}{\partial\tilde{z}}\frac{\partial ln R_0(\tilde{z})}{\partial\tilde{z}}\frac{\partial\tilde{u}_r}{\partial\tilde{r}} \\
	+ \frac{R_0(\tilde{z})^3U_r(\tilde{z})\tilde{r}}{L^3U_z}\frac{\partial ln R_0(\tilde{z})}{\partial\tilde{z}}\frac{\partial ln R_0(\tilde{z})}{\partial\tilde{z}}\frac{\partial\tilde{u}_r}{\partial\tilde{r}}
	- \frac{R_0(\tilde{z})^3U_r(\tilde{z})\tilde{r}}{L^3U_z}\frac{\partial^2 ln R_0(\tilde{z})}{\partial\tilde{z}^2}\frac{\partial\tilde{u}_r}{\partial\tilde{r}} \\
	- \frac{R_0(\tilde{z})^3U_r(\tilde{z})\tilde{r}}{L^3U_z}\frac{\partial ln R_0(\tilde{z})}{\partial\tilde{z}}\frac{\partial^2\tilde{u}_r}{\partial\tilde{r}\partial\tilde{z}} 
	+ \frac{R_0(\tilde{z})^3U_r(\tilde{z})\tilde{r}^2}{L^3U_z}\frac{\partial ln R_0(\tilde{z})}{\partial\tilde{z}}\frac{\partial ln R_0(\tilde{z})}{\partial\tilde{z}}\frac{\partial^2\tilde{u}_r}{\partial\tilde{r}^2} \\
	+ \frac{R_0(\tilde{z})^3}{L^3U_z}\frac{\partial U_r(\tilde{z})}{\partial\tilde{z}}\frac{\partial\tilde{u}_r}{\partial\tilde{z}}
	+ \frac{R_0(\tilde{z})^3U_r(\tilde{z})}{L^3U_z}\frac{\partial^2\tilde{u}_r}{\partial\tilde{z}^2} 
	- \frac{R_0(\tilde{z})^3U_r(\tilde{z})\tilde{r}}{L^3U_z}\frac{\partial ln R_0(\tilde{z})}{\partial\tilde{z}}\frac{\partial^2\tilde{u}_r}{\partial\tilde{z}\partial\tilde{r}} \\
	+ \frac{R_0(\tilde{z})^3\tilde{u}_r}{L^3U_z}\frac{\partial^2 U_r(z)}{\partial\tilde{z}^2} 
	+ \frac{R_0(\tilde{z})^3}{L^3U_z}\frac{\partial U_r(\tilde{z})}{\partial\tilde{z}}\frac{\partial\tilde{u_r}}{\partial\tilde{z}} 
	- \frac{R_0(\tilde{z})^3\tilde{r}}{L^3U_z}\frac{\partial U_r(\tilde{z})}{\partial\tilde{z}}\frac{\partial ln R_0(\tilde{z})}{\partial\tilde{z}}\frac{\partial\tilde{u}_r}{\partial\tilde{r}}\\
	+ \frac{R_0(\tilde{z})U_r(\tilde{z})}{LU_z\tilde{r}} \frac{\partial \tilde{u}_r}{\partial \tilde{r}} 
	- \frac{R_0(\tilde{z})U_r(\tilde{z})\tilde{u}_r}{LU_z\tilde{r}^2}
\Bigg\}.
\end{aligned}
$$
After neglecting all but the first-order $\varepsilon$ terms, the nondimensionalized radial momentum equation simplifies to:
$$
\frac{R_0(\tilde{z}) \tilde{u}_r \tilde{u}_z}{L U_z}\frac{\partial U_r(\tilde{z})}{\partial \tilde{z}}
+ \frac{\partial \tilde{p}}{\partial \tilde{r}}
=0.
$$
Utilizing $U_r(\tilde{z}) = \frac{U_z R_0(\tilde{z})}{L}$, we derive the subsequent equation:
\begin{equation}
	\frac{\tilde{u}_{z} \tilde{u}_{r}}{2 L^{2}} \frac{\partial R_{0}(\tilde{z})^{2}}{\partial\tilde{z}}
	+\frac{\partial \tilde{p}}{\partial \tilde{r}} = 0.
\label{Rad_Moment}
\end{equation}
Note that this equation presents a Darcy-type law, relating the change in pressure and radial velocity.

\subsection{The axial momentum equation in nondimensional form}
Substituting the derived expressions in Sec.~\ref{sec.derivatives} into the Equ.~(\ref{NS}.c), we obtain:
$$
\begin{aligned}
\frac{U_z^2}{L}\frac{\partial\tilde{u}_z}{\partial\tilde{t}} 
+U_r(\tilde{z})\tilde{u}_r\frac{U_z}{R_0(\tilde{z})}\frac{\partial\tilde{u}_z}{\partial\tilde{r}} 
+U_z\tilde{u}_z\left( \frac{U_z}{L}\frac{\partial \tilde{u}_z}{\partial \tilde{z}} 
- \frac{U_z\tilde{r}}{L}\frac{\partial ln R_0(\tilde{z})}{\partial \tilde{z}}\frac{\partial \tilde{u}_z}{\partial \tilde{r}} \right) +\left( \frac{ U_z^2}{L}\frac{\partial\tilde{p}}{\partial\tilde{z}} - \frac{U_z^2\tilde{r}}{L}\frac{\partial lnR_0(\tilde{z})}{\partial\tilde{z}}\frac{\partial\tilde{p}}{\partial\tilde{r}}  \right) \\
= \frac{\mu_f}{\rho_f}\Bigg\{ 
\frac{U_z}{R_0^2(\tilde{z})}\frac{\partial^2\tilde{u}_z}{\partial\tilde{r}^2} 
+ \frac{U_z}{L^2}\frac{\partial^2\tilde{u}_z}{\partial\tilde{z}^2} - \frac{U_z\tilde{r}}{L^2}\frac{\partial ln R_0(\tilde{z})}{\partial\tilde{z}}\frac{\partial^2\tilde{u}_z}{\partial\tilde{z}\partial\tilde{r}} 
+ \frac{U_z\tilde{r}}{L^2}\frac{\partial ln R_0(\tilde{z})}{\partial\tilde{z}}\frac{\partial ln R_0(\tilde{z})}{\partial\tilde{z}}\frac{\partial\tilde{u}_z}{\partial\tilde{r}}\\
- \frac{U_z\tilde{r}}{L^2}\frac{\partial^2ln R_0(\tilde{z})}{\partial\tilde{z}^2}\frac{\partial\tilde{u}_z}{\partial\tilde{r}}
+ \frac{U_z\tilde{r}^2}{L^2}\frac{\partial ln R_0(\tilde{z})}{\partial\tilde{z}}\frac{\partial ln R_0(\tilde{z})}{\partial\tilde{z}}\frac{\partial^2\tilde{u}_z}{\partial\tilde{r}^2} - \frac{U_z\tilde{r}}{L^2}\frac{\partial^2\tilde{u}_z}{\partial\tilde{r}\partial\tilde{z}} 
+ \frac{1}{R_0(\tilde{z})\tilde{r}}\frac{U_z}{R_0(\tilde{z})}\frac{\partial\tilde{u}_z}{\partial\tilde{r}} \Bigg\}.
\end{aligned}
$$
Dividing by $U_{z}^{2}$ and then multiplying by $L$ yields:
$$
\begin{aligned}
\frac{\partial\tilde{u}_z}{\partial\tilde{t}} 
+\frac{U_r(\tilde{z})L\tilde{u_r}}{R_0(\tilde{z})U_z}\frac{\partial\tilde{u}_z}{\partial\tilde{r}} 
+\tilde{u}_z\frac{\partial \tilde{u}_z}{\partial \tilde{z}} 
- \tilde{u}_z\tilde{r}\frac{\partial ln R_0(\tilde{z})}{\partial \tilde{z}}\frac{\partial \tilde{u}_z}{\partial \tilde{r}}  
+\frac{\partial\tilde{p}}{\partial\tilde{z}} 
- \frac{\partial ln R_0(\tilde{z})}{\partial\tilde{z}}\frac{\partial\tilde{p}}{\partial\tilde{r}} \\
= \frac{\mu_f L}{\rho_f R_0(\tilde{z})^2U_z}\Bigg\{ 
\frac{\partial^2\tilde{u}_z}{\partial\tilde{r}^2} 
+ \frac{R_0(\tilde{z})^2}{L^2}\frac{\partial^2\tilde{u}_z}{\partial\tilde{z}^2} 
- \frac{R_0(\tilde{z})^2\tilde{r}}{L^2}\frac{\partial ln R_0(\tilde{z})}{\partial\tilde{z}}\frac{\partial^2\tilde{u}_z}{\partial\tilde{z}\partial\tilde{r}} 
+ \frac{R_0(\tilde{z})^2\tilde{r}}{L^2}\frac{\partial ln R_0(\tilde{z})}{\partial\tilde{z}}\frac{\partial ln R_0(\tilde{z})}{\partial\tilde{z}}\frac{\partial\tilde{u}_z}{\partial\tilde{r}} \\
- \frac{R_0(\tilde{z})^2\tilde{r}}{L^2}\frac{\partial^2ln(R_0(z))}{\partial\tilde{z}^2}\frac{\partial\tilde{u}_z}{\partial\tilde{r}}
+ \frac{R_0(\tilde{z})^2\tilde{r}^2}{L^2}\frac{\partial ln R_0(\tilde{z})}{\partial\tilde{z}}\frac{\partial ln R_0(\tilde{z})}{\partial\tilde{z}}\frac{\partial^2\tilde{u}_z}{\partial\tilde{r}^2} 
- \frac{R_0(\tilde{z})^2\tilde{r}}{L^2}\frac{\partial^2\tilde{u}_z}{\partial\tilde{r}\partial\tilde{z}} 
+ \frac{1}{\tilde{r}}\frac{\partial\tilde{u}_z}{\partial\tilde{r}} 
\Bigg\}.
\end{aligned}
$$
Noticing that several terms are of order $\varepsilon^{2}$, we can neglect them, which simplifies the equation to:
$$
\begin{aligned}
\frac{\partial\tilde{u}_z}{\partial\tilde{t}} 
+\tilde{u}_r\frac{\partial\tilde{u}_z}{\partial\tilde{r}} 
+\tilde{u}_z\frac{\partial \tilde{u}_z}{\partial \tilde{z}} 
- \tilde{u}_z\tilde{r}\frac{\partial ln R_0(\tilde{z})}{\partial \tilde{z}}\frac{\partial \tilde{u}_z}{\partial \tilde{r}}  
+\frac{\partial\tilde{p}}{\partial\tilde{z}} - \frac{\partial ln R_0(\tilde{z})}{\partial\tilde{z}}\frac{\partial\tilde{p}}{\partial\tilde{r}}
= \frac{\nu_fL}{R_0(\tilde{z})^2U_z}\frac{1}{\tilde{r}}\frac{\partial}{\partial\tilde{r}}\left(\tilde{r}\frac{\partial\tilde{u}_z}{\partial\tilde{z}} \right)
\end{aligned}
$$
We define the Reynolds number as:
$$
\mathrm{Re}=\frac{\rho_f U_{z} R_{0}^{2}}{\mu_f L},
$$
so that the leading-order axial momentum equation in non-dimensional form becomes:
\begin{equation}
\begin{aligned}
\frac{\partial\tilde{u}_z}{\partial\tilde{t}} 
+\tilde{u}_r\frac{\partial\tilde{u}_z}{\partial\tilde{r}} 
+\tilde{u}_z\frac{\partial \tilde{u}_z}{\partial \tilde{z}} 
- \tilde{u}_z\tilde{r}\frac{\partial ln R_0(\tilde{z})}{\partial \tilde{z}}\frac{\partial \tilde{u}_z}{\partial \tilde{r}}  
+\frac{\partial\tilde{p}}{\partial\tilde{z}} 
- \frac{\partial ln R_0(\tilde{z})}{\partial\tilde{z}}\frac{\partial\tilde{p}}{\partial\tilde{r}}
= \frac{1}{Re}\frac{1}{\tilde{r}}\frac{\partial}{\partial\tilde{r}}\left(\tilde{r}\frac{\partial\tilde{u}_z}{\partial\tilde{z}} \right)
\end{aligned}
\label{Axial1}
\end{equation}
Noticing that:
$$
\frac{\partial}{\partial \tilde{r}} (\tilde{r} \tilde{u}_{r} \tilde{u}_{z}) 
 =\tilde{u}_{z} \frac{\partial}{\partial \tilde{r}} (\tilde{r} \tilde{u}_{r})
+\tilde{r} \tilde{u}_{r} \frac{\partial \tilde{u}_{z}}{\partial \tilde{r}}, \quad
\frac{\partial}{\partial \tilde{z}} (\tilde{r} \tilde{u}_{z}^{2}) =\tilde{u}_{z} \frac{\partial}{\partial \tilde{z}} (\tilde{r} \tilde{u}_{z})
+\tilde{r} \tilde{u}_{z} \frac{\partial \tilde{u}_{z}}{\partial \tilde{z}},
$$
we utilize Equ.~(\ref{mass_conser}) to obtain:
 $$
\frac{\partial}{\partial \tilde{r}} (\tilde{r} \tilde{u}_{r} \tilde{u}_{z})
+ \frac{\partial}{\partial \tilde{z}} (\tilde{r} \tilde{u}_{z}^{2})
= \tilde{r} \tilde{u}_{r} \frac{\partial \tilde{u}_{z}}{\partial \tilde{r}}
+ \tilde{r} \tilde{u}_{z} \frac{\partial \tilde{u}_z}{\partial \tilde{z}}
+ \tilde{u}_z \tilde{r}^2\frac{\partial ln R_0(\tilde{z})}{\partial\tilde{z}}\frac{\partial\tilde{u}_z}{\partial\tilde{r}}.
 $$
By multiplying Equ.~(\ref{Axial1}) by $\tilde{r}$ and then applying the relationships stated earlier along with Equ.~(\ref{mass_conser}), we obtain the axial momentum equation in conservation form:
\begin{equation}
\begin{split}
\frac{\partial }{\partial \tilde{t}}(\tilde{r} \tilde{u}_{z})
		+\frac{\partial}{\partial \tilde{r}} (\tilde{r} \tilde{u}_{r} \tilde{u}_{z})
		+\frac{\partial}{\partial \tilde{z}} (\tilde{r} \tilde{u}_{z}^{2})
		- 2\tilde{u}_z\tilde{r}^2\frac{\partial ln R_0(\tilde{z})}{\partial\tilde{z}}\frac{\partial\tilde{u}_z}{\partial\tilde{r}}
		+\frac{\partial}{\partial \tilde{z}} (\tilde{r} \tilde{p})
		- \tilde{r}\frac{\partial ln R_0(\tilde{z})}{\partial\tilde{z}}\frac{\partial\tilde{p}}{\partial\tilde{r}}\\
 =\frac{1}{\operatorname{Re}} \frac{\partial}{\partial \tilde{r}}\left(\tilde{r} \frac{\partial \tilde{u}_{z}}{\partial \tilde{r}}\right).
\label{axial_moment}
\end{split}
\end{equation}

\subsection{The leading order system and its averaged equations}
In summary, we derive the following leading-order system:
\begin{equation}
\left\{
    \begin{aligned}
        &\frac{\partial }{\partial \tilde{r}}(\tilde{r} \tilde{u}_{r})
         + \frac{\partial }{\partial \tilde{z}} (\tilde{r} \tilde{u}_{z})
        - \tilde{r}^2\frac{\partial \ln R_{0}(\tilde{z})}{\partial \tilde{z}}\frac{\partial \tilde{u}_z}{\partial \tilde{r}}
        = 0, \\
        &\frac{\partial }{\partial \tilde{t}}(\tilde{r} \tilde{u}_{z})
        +\frac{\partial}{\partial \tilde{r}} (\tilde{r} \tilde{u}_{r} \tilde{u}_{z})
        +\frac{\partial}{\partial \tilde{z}} (\tilde{r} \tilde{u}_{z}^{2})
        - 2\tilde{u_z}\tilde{r}^2\frac{\partial ln R_0(\tilde{z})}{\partial\tilde{z}}\frac{\partial\tilde{u}_z}{\partial\tilde{r}}
        +\frac{\partial}{\partial \tilde{z}} (\tilde{r} \tilde{p}) 
        - \tilde{r}\frac{\partial ln R_0(\tilde{z})}{\partial\tilde{z}}\frac{\partial\tilde{p}}{\partial\tilde{r}}\\
       & \hspace{4in}=\frac{1}{\operatorname{Re}} \frac{\partial}{\partial \tilde{r}}\left(\tilde{r} \frac{\partial \tilde{u}_{z}}{\partial \tilde{r}}\right), \\
        &\frac{\tilde{u}_{z} \tilde{u}_{r}}{2 L^{2}} \frac{\partial R_{0}^{2}(\tilde{z})}{\partial\tilde{z}}
        +\frac{\partial \tilde{p}}{\partial \tilde{r}} = 0.
\end{aligned}
\right.
\label{1D_system_1}
\end{equation}
To obtain the averaged equation, we integrate the leading-order system (\ref{1D_system_1}) across the cross-section of the blood vessel. This integration is performed with respect to the radius 
$\tilde{r}$ from $\tilde{r}=0$ to $\tilde{r}=\tilde{R}(\tilde{z}, \tilde{t})$. Below, we will derive the averaged equations step by step.
	
\subsubsection{The conservation of mass in reduced form}
First, we address the Equ.~(\ref{1D_system_1}.a) by integrating with respect to the radius, we obtain:
\begin{equation}
\int_{0}^{\tilde{R}} \frac{\partial }{\partial \tilde{r}} (\tilde{r}\tilde{u}_{r}) d \tilde{r}
+\int_{0}^{\tilde{R}} \frac{\partial}{\partial \tilde{z}} (\tilde{r} \tilde{u}_{z}) d \tilde{r}
-\int_{0}^{\tilde{R}} \tilde{r}^2\frac{\partial \ln R_{0}(\tilde{z})}{\partial \tilde{z}}\frac{\partial \tilde{u}_z}{\partial \tilde{r}} d \tilde{r}
=0. 
\label{com_integral}
\end{equation}
The second term in the equation above can be rewritten using the Leibniz integral rule as follows:
$$
\int_{0}^{\tilde{R}(\tilde{z}, \tilde{t})} \frac{\partial}{\partial \tilde{z}} (\tilde{r} \tilde{u}_{z}) d \tilde{r}
= \frac{\partial}{\partial \tilde{z}} \int_{0}^{\tilde{R}(\tilde{z}, \tilde{t})} (\tilde{r} \tilde{u}_{z}) d \tilde{r}
-\frac{\partial \tilde{R}(\tilde{z}, \tilde{t})}{\partial \tilde{z}}
\tilde{R}(\tilde{z}, \tilde{t}) \tilde{u}_z (\tilde{R}, \tilde{z}, \tilde{t}).
$$
Concerning the third term, we have:
$$
\begin{aligned}
-\int_{0}^{\tilde{R}} \tilde{r}^2\frac{\partial \ln R_{0}(\tilde{z})}{\partial \tilde{z}}\frac{\partial \tilde{u}_z}{\partial \tilde{r}} d \tilde{r}
&= -\frac{\partial \ln R_{0}(\tilde{z})}{\partial \tilde{z}}\left(
\tilde{r}^2\tilde{u}_z\Bigg|_0^{\tilde{R}} 
- \int_0^{\tilde{R}}2\tilde{r}\tilde{u}_zd\tilde{r}
\right) \\
&= - \frac{\partial ln R_0(\tilde{z})}{\partial\tilde{z}}\tilde{R}^2(\tilde{z}, \tilde{t})\tilde{u}_z(\tilde{R}, \tilde{z}, \tilde{t})
+ 2\frac{\partial ln R_0(\tilde{z})}{\partial\tilde{z}}\int_0^{\tilde{R}}\tilde{r}\tilde{u}_zd\tilde{r}
\end{aligned}
$$ 
Therefore, Equ.~(\ref{com_integral}) becomes:
\begin{equation}
\begin{split}
\begin{aligned}
\tilde{R}(\tilde{z}, \tilde{t}) \tilde{u}_{r}(\tilde{R}, \tilde{z}, \tilde{t})
+\frac{\partial}{\partial \tilde{z}} \int_{0}^{\tilde{R}(\tilde{z}, \tilde{t})} \tilde{r}\tilde{u}_{z}  d \tilde{r}
-\frac{\partial \tilde{R}(\tilde{z}, \tilde{t})}{\partial \tilde{z}} \tilde{R}(\tilde{z}, \tilde{t}){\tilde{u}_{z}(\tilde{R}, \tilde{z}, \tilde{t})}\\
- \frac{\partial ln R_0(\tilde{z})}{\partial\tilde{z}}\tilde{R}^2(\tilde{z}, \tilde{t})\tilde{u}_z(\tilde{R}, \tilde{z}, \tilde{t})
+ 2\frac{\partial ln R_0(\tilde{z})}{\partial\tilde{z}}\int_0^{\tilde{R}}\tilde{r}\tilde{u}_zd\tilde{r}
=0.
\end{aligned}
\label{mass_conser_int}
\end{split}
\end{equation}
Recall the no-slip condition at the interface ${\Gamma(t)}=(\tilde{R}, \tilde{z}, \tilde{t})$, satisfying that:
	$$
	\left.\left(u_{r}, u_{z}\right)\right|_{\Gamma(t)}=\left(\frac{\partial \eta_{r}}{\partial t}, 0\right)=\left(\frac{\partial R}{\partial t}, 0\right),
	$$
Thus, in non-dimensional form, we state:
	$$
	\left.\left(\tilde{u}_{r}, \tilde{u}_{z}\right)\right|_{\Gamma(t)}=\left(\frac{\partial \tilde{R}}{\partial \tilde{t}}, 0\right).
	$$
Take this into account, Equ.~(\ref{mass_conser_int}) becomes:
$$
\tilde{R}(\tilde{z}, \tilde{t}) \frac{\partial \tilde{R}(\tilde{z}, \tilde{t})}{\partial \tilde{t}}
+\frac{\partial}{\partial \tilde{z}} \int_{0}^{\tilde{R}(\tilde{z}, \tilde{t})}\tilde{r} \tilde{u}_{z}  d \tilde{r}
+2\frac{\partial \ln R_{0}(\tilde{z})}{\partial \tilde{z}} \int_{0}^{\tilde{R}(\tilde{z}, \tilde{t})} \tilde{r} \tilde{u}_{z} d \tilde{r} 
=0.
$$
Next, we introduce notation for the average axial component of velocity over the cross-section as:
\begin{equation}
\tilde{U}=\frac{2}{\tilde{R}^{2}} \int_{0}^{\tilde{R}}  \tilde{r}\tilde{u}_{z} d \tilde{r},
\label{axial_aver_vel_nondimensional}
\end{equation}
so we have:
$$
\frac{\partial \tilde{R}^{2}}{\partial \tilde{t}}
+\frac{\partial}{\partial \tilde{z}} (\tilde{R}^{2} \tilde{U})
+2\frac{\partial \ln R_{0}(\tilde{z})}{\partial \tilde{z}} ({\tilde{R}^{2}\tilde{U}}) 
=0.
$$
Introducing the following equation for the scaled cross-sectional area $\tilde{A}=\tilde{R}^{2}$ and flow rate $\tilde{Q}=\tilde{A} \tilde{U}$,
we obtain:
\begin{equation}
\begin{split}
\frac{\partial \tilde{A}}{\partial \tilde{t}}
+\frac{\partial\tilde{Q}}{\partial \tilde{z}}
+2\frac{\partial \ln R_{0}(\tilde{z})}{\partial \tilde{z}} {\tilde{Q}} 
=0.
\end{split}
\label{reduced_conser_mass}
\end{equation}
	
\subsubsection{The axial momentum equation in reduced form}
We address Equ.~(\ref{1D_system_1}.b) by integrating with respect to the radius to obtain:
\begin{equation}
    \begin{aligned}
        \int_{0}^{\tilde{R}} \frac{\partial }{\partial \tilde{t}}(\tilde{r} \tilde{u}_{z}) d \tilde{r}
        +\int_{0}^{\tilde{R}} \frac{\partial}{\partial \tilde{r}} (\tilde{r} \tilde{u}_{r} \tilde{u}_{z}) d \tilde{r}
        +\int_{0}^{\tilde{R}} \frac{\partial}{\partial \tilde{z}} (\tilde{r} \tilde{u}_{z}^{2}) d \tilde{r}
        - \int_{0}^{\tilde{R}} 2\tilde{u}_z\tilde{r}^2\frac{\partial ln R_0(\tilde{z})}{\partial\tilde{z}}\frac{\partial\tilde{u}_z}{\partial\tilde{r}} d \tilde{r}\\
        +\int_{0}^{\tilde{R}} \frac{\partial}{\partial \tilde{z}} (\tilde{r} \tilde{p}) d \tilde{r} 
        - \int_{0}^{\tilde{R}} \tilde{r}\frac{\partial ln R_0(\tilde{z})}{\partial\tilde{z}}\frac{\partial\tilde{p}}{\partial\tilde{r}} d \tilde{r}
         =\int_{0}^{\tilde{R}} \frac{1}{\operatorname{Re}} \frac{\partial}{\partial \tilde{r}}\left(\tilde{r} \frac{\partial \tilde{u}_{z}}{\partial \tilde{r}}\right) d \tilde{r}.
    \end{aligned}
    \label{axial_momentum_reduced_2}
\end{equation}
For the fourth term in Equ.~(\ref{axial_momentum_reduced_2}), we observe:
$$
\begin{aligned}
- \int_{0}^{\tilde{R}} 2\tilde{u}_z\tilde{r}^2\frac{\partial lnR_0(\tilde{z})}{\partial\tilde{z}}\frac{\partial\tilde{u}_z}{\partial\tilde{r}} d \tilde{r}
 &= -\frac{\partial ln R_0(\tilde{z})}{\partial\tilde{z}}\int_0^{\tilde{R}} 2\tilde{u}_z\tilde{r}^2\frac{\partial\tilde{u}_z}{\partial\tilde{r}} d \tilde{r} \\
&= -\frac{\partial ln R_0(\tilde{z})}{\partial\tilde{z}}
\left(
2\tilde{u}_z^2\tilde{r}^2\Bigg|_0^{\tilde{R}}
- 4\int_0^{\tilde{R}}\tilde{r}\tilde{u}_z^2 d \tilde{r}
-2\int_0^{\tilde{R}}\tilde{r}^2\tilde{u}_z\frac{\partial\tilde{u}_z}{\partial\tilde{r}}d \tilde{r}
	 \right).
\end{aligned}
$$
We observe a repetition term on the right-hand side of the equation above. By moving it to the left-hand side, we obtain:
$$
\begin{aligned}
	-\frac{\partial ln R_0(\tilde{z})}{\partial\tilde{z}}\int_0^{\tilde{R}} \tilde{u}_z\tilde{r}^2\frac{\partial\tilde{u}_z}{\partial\tilde{r}} d \tilde{r}
	= -\frac{1}{2}\frac{\partial ln R_0(\tilde{z})}{\partial\tilde{z}}
 \tilde{R}^2(\tilde{z}, \tilde{t})\tilde{u}_z(\tilde{R}, \tilde{z}, \tilde{t})
	+ \frac{\partial ln(R_0(z))}{\partial\tilde{z}}\int_0^{\tilde{R}} \tilde{r}\tilde{u}_z^2 d\tilde{r}.
\end{aligned}
$$
\color{black}
To handle the sixth term in Equ.~(\ref{axial_momentum_reduced_2}), we introduce a variable name $I$ to represent the integral below:
$$
\mathcal{I}=-\int_0^{\tilde{R}} \tilde{r} \frac{\partial \ln R_0(\tilde{z})}{\partial \tilde{z}} \frac{\partial \tilde{p}}{\partial \tilde{r}}d\tilde{r}
=-\frac{\partial \ln R_0(\tilde{z})}{\partial \tilde{z}} \int_0^{\tilde{R}} \tilde{r} \frac{\partial \tilde{p}}{\partial \tilde{r}} d\tilde{r}
$$
Using Equation (\ref{1D_system_1}.c), we derive:
\begin{equation}
\mathcal{I}= \frac{\partial \ln R_0(\tilde{z})}{\partial \tilde{z}} 
\frac{1}{2 L^2} \frac{\partial R_0^2(\tilde{z})}{\partial \tilde{z}}
\int_0^{\tilde{R}} \tilde{r} \tilde{u}_z \tilde{u}_r d \tilde{r}.
\label{Intergral_6}
\end{equation}
To compute the integral $\int_0^{\tilde{R}} \tilde{r} \tilde{u}_z \tilde{u}_r d \tilde{r}$, we will achieve it in two steps. The first step involves finding the expression of $\tilde{r} \tilde{u}_r$. By assuming a parabolic velocity profile of \(\tilde{u}_z = 2 \tilde{U} \left(1 - \frac{\tilde{r}^2}{\tilde{R}^2}\right)\), which depends solely on \(\tilde{r}\), and from Equation (\ref{1D_system_1}.a), we know that:
$$
\frac{\partial \tilde{r} \tilde{u}_r}{\partial \tilde{r}}= -\frac{\partial \tilde{r} \tilde{u}_z}{\partial \tilde{z}} 
+\tilde{r}^2 \frac{\partial ln R_0(\tilde{z})}{\partial \tilde{z}} \frac{\partial \tilde{u}_z}{\partial \tilde{r}}.
$$
Integrate above equation with respect to $\tilde{r}$ from 0 to $\tilde{r}$, yields:
$$
\begin{aligned}
\tilde{r} \tilde{u}_r=\int_0^{\tilde{r}}\frac{\partial \tilde{r} \tilde{u}_r}{\partial \tilde{r}} d\tilde{r}
=&-\int_0^{\tilde{r}}\frac{\partial \tilde{r} \tilde{u}_z}{\partial \tilde{z}} d\tilde{r} 
+\int_0^{\tilde{r}} \tilde{r}^2 \frac{\partial ln R_0(\tilde{z})}{\partial \tilde{z}} \frac{\partial \tilde{u}_z}{\partial \tilde{r}} d\tilde{r}\\
=&-\frac{\partial}{\partial \tilde{z}}\int_0^{\tilde{r}}\tilde{r} \tilde{u}_z d \tilde{r}
+\frac{\partial \ln R_0(\tilde{z})}{\partial z} \int_0^{\tilde{r}} \tilde{r}^2 \frac{\partial\left[2\tilde{U}\left(1-\frac{\tilde{r}^2}{\tilde{R}^2}\right)\right] }{\partial \tilde{r}} d \tilde{r}\\
=&-\frac{\partial}{\partial \tilde{z}} \int_0^{\tilde{r}} 2 \tilde{r}\tilde{U} \left(1-\frac{\tilde{r}^2}{\tilde{R}^2}\right) d \tilde{r}
+\frac{\partial \ln R_0(\tilde{z})}{\partial \tilde{z}} \int_0^{\tilde{r}} 2 \tilde{U} \tilde{r}^2\left(\frac{-2 \tilde{r}}{\tilde{R}^2}\right) d \tilde{r} \\
=&-\frac{\partial}{\partial \tilde{z}}\left[2 \tilde{U}\left(\frac{\tilde{r}^2}{2}-\frac{\tilde{r}^4}{4 \tilde{R}^2}\right)\right]+\frac{\partial \ln R_0(\tilde{z})}{\partial \tilde{z}}\left[2 \tilde{U}\left(-\frac{\tilde{r}^4}{2 \tilde{R}^2}\right)\right] \\
=&-2 \frac{\partial \tilde{U}}{\partial \tilde{z}}\left(\frac{\tilde{r}^2}{2}-\frac{\tilde{r}^4}{4 \tilde{R}^2}\right)
+2 \tilde{U}\left(\frac{\partial\left(\frac{\tilde{r}^4}{4 \tilde{R}^2}\right)}{\partial \tilde{z}}\right)
+\frac{\partial \ln R_0(\tilde{z})}{\partial \tilde{z}} 2 \tilde{U}\left(-\frac{\tilde{r}^4}{2 \tilde{R}^2}\right) \\
=&-\left(\tilde{r}^2-\frac{\tilde{r}^4}{2 \tilde{R}^2}\right) \frac{\partial \tilde{U}}{\partial \tilde{z}}
-\frac{\tilde{U}\tilde{r}^4}{\tilde{R}^3}\frac{\partial \tilde{R}}{\partial \tilde{z}}
-\frac{\partial \ln R_0(\tilde{z})}{\partial \tilde{z}}
\frac{\tilde{U} \tilde{r}^4}{\tilde{R}^2} \\
=&-\tilde{r}^2\left(1-\frac{\tilde{r}^2}{2 \tilde{R}^2}\right) \frac{\partial \tilde{U}}{\partial \tilde{z}}
-\frac{\tilde{U} \tilde{r}^4}{\tilde{R}^2} \left(\frac{\partial \ln \tilde{R}}{\partial \tilde{z}}
+\frac{\partial \ln R_0(\tilde{z})}{\partial \tilde{z}} \right) 
\end{aligned}
$$
Therefore,
$$
\begin{aligned}
\tilde{r} \tilde{u}_r \tilde{u}_z=2 \tilde{U}\left(1-\frac{\tilde{r}^2}{\tilde{R}^2}\right)\left(-\tilde{r}^2\left(1-\frac{\tilde{r}^2}{2 \tilde{R}^2}\right) \frac{\partial \tilde{U}}{\partial \tilde{z}}-\frac{\tilde{U} \tilde{r}^4}{\tilde{R}^2} \left( \frac{\partial \ln \tilde{R}}{\partial \tilde{z}}
+\frac{\partial \ln R_0(\tilde{z})}{\partial \tilde{z}} \right) \right) .
\end{aligned}
$$
Note that terms such as $\tilde{U}$, 
$\frac{\partial \tilde{U}}{\partial \tilde{z}}$, $\frac{\partial \ln R_0(\tilde{z})}{\partial \tilde{z}}$ and $\frac{\partial \ln \tilde{R}}{\partial \tilde{z}}$ do not depend on $\tilde{r}$. Therefore, we can pull them out of the integration, which yields:
\begin{equation}
\begin{aligned}
 &\int_0^{\tilde{R}} \tilde{r} \tilde{u}_r \tilde{u}_z d \tilde{r}\\
&=-2 \tilde{U} \frac{\partial \tilde{U}}{\partial \tilde{z}} \int_0^{\tilde{R}}
\left(1-\frac{\tilde{r}^2}{\tilde{R}^2}\right) \tilde{r}^2\left(1-\frac{\tilde{r}^2}{2 \tilde{R}^2}\right) d \tilde{r} 
- 2\tilde{U}^2 \left(\frac{\partial \ln \tilde{R}}{\partial \tilde{z}}
+\frac{\partial \ln R_0(\tilde{z})}{\partial \tilde{z}} \right) \int_0^{\tilde{R}}\left(1-\frac{\tilde{r}^2}{\tilde{R}^2}\right) \frac{\tilde{r}^4}{\tilde{R}^2} d \tilde{r} \\
& =-2 \tilde{U} \frac{\partial \tilde{U}}{\partial \tilde{z}}\left[
\frac{\frac{\tilde{r}^7}{14}-\frac{3 \tilde{R}^2 \tilde{r}^5}{10}}{\tilde{R}^4}+\frac{\tilde{r}^3}{3}
\right]_{\tilde{r}=0}^{\tilde{r}=\tilde{R}} 
-2 \tilde{U}^2 \left(\frac{\partial \ln \tilde{R}}{\partial \tilde{z}}
+\frac{\partial \ln R_0(\tilde{z})}{\partial \tilde{z}} \right)\left[\frac{\tilde{r}^5}{5 \tilde{R}^2}-\frac{\tilde{r}^7}{7 \tilde{R}^4}\right]_{\tilde{r}=0}^{\tilde{r}=\tilde{R}} \\
& =-2 \tilde{U} \frac{\partial \tilde{U}}{\partial \tilde{z}}\left(\frac{11}{105} \tilde{R}^3\right)-2 \tilde{U}^2 \left(\frac{\partial \ln \tilde{R}}{\partial \tilde{z}}
+\frac{\partial \ln R_0(\tilde{z})}{\partial \tilde{z}} \right)\left(\frac{2}{35} \tilde{R}^3\right) \\
& \approx -(\frac{2}{35}\tilde{R}^3)\left[
\frac{\partial \tilde{U}^2}{\partial \tilde{z}}
+2\tilde{U}^2\left(\frac{\partial ln R_0(\tilde{z})}{\partial \tilde{z}}
+\frac{\partial ln \tilde{R}}{\partial \tilde{z}}\right)
\right].
\end{aligned}
\label{approx}
\end{equation}
Therefore, the integration $\mathcal{I}$ as shown in Equ.~(\ref{Intergral_6}) becomes:
\begin{equation}
\begin{aligned}
\mathcal{I}=&-\frac{\partial \ln R_0(\tilde{z})}{\partial \tilde{z}} 
\frac{1}{2 L^2} \frac{\partial R_0^2(\tilde{z})}{\partial \tilde{z}}
(\frac{2}{35}\tilde{R}^3)
\left[
\frac{\partial \tilde{U}^2}{\partial \tilde{z}}
+2\tilde{U}^2\left(\frac{\partial ln R_0(\tilde{z})}{\partial \tilde{z}}
+\frac{\partial ln \tilde{R}}{\partial \tilde{z}}\right)
\right]\\
=&-\frac{2}{35L^2}\left(\frac{\partial R_0(\tilde{z})}{\partial \tilde{z}}\right)^2
\frac{\tilde{A}^2}{\tilde{R}}
\left[
\frac{\partial \tilde{U}^2}{\partial \tilde{z}}
+2\tilde{U}^2\left(\frac{\partial ln R_0(\tilde{z})}{\partial \tilde{z}}
+2\frac{\partial ln \tilde{R}}{\partial \tilde{z}}\right)
\right]\\
=&-\frac{2}{35L^2}\left(\frac{\partial R_0(\tilde{z})}{\partial \tilde{z}}\right)^2
\frac{1}{\tilde{R}}
\left[
\frac{\partial \tilde{Q}^2}{\partial \tilde{z}}
-\tilde{U}^2\frac{\partial \tilde{A}^2}{\partial \tilde{z}}
+2\tilde{Q}^2\left(\frac{\partial ln R_0(\tilde{z})}{\partial \tilde{z}}
+\frac{\partial ln \tilde{R}}{\partial \tilde{z}}\right)
\right]\\ 
=&-\frac{2}{35L^2}\left(\frac{\partial R_0(\tilde{z})}{\partial \tilde{z}}\right)^2
\frac{1}{\tilde{R}}
\left[
\frac{\partial \tilde{Q}^2}{\partial \tilde{z}}
-\tilde{U}^2\frac{\partial \tilde{A}^2}{\partial \tilde{z}}
+2\tilde{Q}^2\frac{\partial ln (R_0(\tilde{z})\tilde{R})}{\partial \tilde{z}}
\right]\\
\end{aligned}
\end{equation}
We emphasize that the term $\mathcal{I}$ can be further simplified when converted to its dimensional form, see details in Sec.~\ref{retrieve_dimensional_form}. 
For the remainder of this section, we will leave it in its nondimensional form for now.

\color{black}

We pull out derivatives from the integral where possible, bearing in mind that $\tilde{R}=\tilde{R}(\tilde{z}, \tilde{t})$ and $\tilde{p}$ is independent of $\tilde{r}$, we deduce:
\begin{equation}
	\begin{aligned}
\frac{\partial}{\partial \tilde{t}} \int_{0}^{\tilde{R}} \tilde{r} \tilde{u}_{z} d \tilde{r}
-\frac{\partial \tilde{R}}{\partial \tilde{t}} \tilde{R} {\tilde{u}_{z}(\tilde{R}, \tilde{z}, \tilde{t})}
+\tilde{R} \tilde{u}_{r}(\tilde{R}, \tilde{z}, \tilde{t}) {\tilde{u}_{z}(\tilde{R}, \tilde{z}, \tilde{t})}
+\frac{\partial}{\partial \tilde{z}} \int_{0}^{\tilde{R}} \tilde{r} \tilde{u}_{z}^{2} d \tilde{r}
-\frac{\partial \tilde{R}}{\partial \tilde{z}} \tilde{R} {\tilde{u}_{z}^{2}(\tilde{R}, \tilde{z}, \tilde{t})}\\
+2\frac{\partial ln R_0(\tilde{z})}{\partial\tilde{z}}\int_0^{\tilde{R}} \tilde{r}\tilde{u}_z^2 d\tilde{r}
- \frac{\partial ln R_0(\tilde{z})}{\partial\tilde{z}}\tilde{R}^2\tilde{u}_z(\tilde{R}, \tilde{z}, \tilde{t})
+\frac{\partial\tilde{p}}{\partial\tilde{z}}\frac{\tilde{R}^2}{2}
\textcolor{black}{+\mathcal{I}}
=\frac{1}{\mathrm{Re}} \tilde{R} \frac{\partial \tilde{u}_{z}(\tilde{R}, \tilde{z}, \tilde{t})}{\partial \tilde{r}}.
	\end{aligned}
\end{equation}
Noticing that the terms involving $\tilde{u}_z(\tilde{R}, \tilde{z}, \tilde{t})$ are zero, Equation (\ref{1D_system_1}.b) can be expressed as follows:
\begin{equation}
    \frac{\partial}{\partial \tilde{t}} \int_{0}^{\tilde{R}} \tilde{r} \tilde{u}_{z} d \tilde{r}
    +\frac{\partial}{\partial \tilde{z}} \int_{0}^{\tilde{R}} \tilde{r} \tilde{u}_{z}^{2} d \tilde{r}
	+ 2\frac{\partial ln(R_0(z))}{\partial\tilde{z}}\int_0^{\tilde{R}} \tilde{r}\tilde{u}_z^2 d\tilde{r}
    +\frac{\tilde{R}^2}{2}\frac{\partial\tilde{p}}{\partial\tilde{z}}
    \textcolor{black}{+\mathcal{I}}
    =\frac{1}{\mathrm{Re}} \tilde{R} \frac{\partial \tilde{u}_{z}(\tilde{R}, \tilde{z}, \tilde{t})}{\partial \tilde{r}}.
\end{equation}
Next, we multiply the equation above by 2, and introduce the term:
$$
\alpha=\frac{2}{\tilde{R}^{2} \tilde{U}^{2}} \int_{0}^{\tilde{R}} \tilde{r} \tilde{u}_{z}^{2} d \tilde{r},
$$
which represents the quadratic velocity term arising from nonlinear advection, commonly referred to as the correction factor or the Coriolis coefficient. Notably, for the Poiseuille velocity profile, $\alpha$ remains constant.

The averaged axial momentum equation can be rewritten in terms of  $\alpha$, $\tilde{U}$, $\tilde{A}$, and $\tilde{Q}$:
\begin{equation}
\begin{split}
\frac{\partial \tilde{Q}}{\partial \tilde{t}}
+\frac{\partial}{\partial \tilde{z}} \left(\alpha \frac{\tilde{Q}^{2}}{\tilde{A}}\right)
+ 2\alpha\frac{\tilde{Q}^2}{\tilde{A}}\frac{\partial ln R_0(\tilde{z})}{\partial\tilde{z}}
+\tilde{A}\frac{\partial\tilde{p}}{\partial\tilde{z}}
\textcolor{black}{+\mathcal{I}}
=\frac{2}{\operatorname{Re}}\tilde{R} \frac{\partial \tilde{u}_{z}(\tilde{R}, \tilde{z}, \tilde{t})}{\partial \tilde{r}}.
\end{split}
\label{eq18}
\end{equation}

\subsection{Reduced $(\tilde{A}, \tilde{Q})$ system and its reconstruction to the original dimensional form}
\label{retrieve_dimensional_form}
By combining Equs.~(\ref{reduced_conser_mass}) and (\ref{eq18}), we obtain the following non-dimensional equations based on key variables ($\tilde{A}$ and $\tilde{Q}$):
\begin{equation}
\left\{
	\begin{aligned}
	&\frac{\partial \tilde{A}}{\partial \tilde{t}}
	+\frac{\partial \tilde{Q}}{\partial \tilde{z}}
	+2\tilde{Q} \frac{\partial \ln R_{0}(\tilde{z})}{\partial \tilde{z}}
	=0, \\
	&\frac{\partial \tilde{Q}}{\partial \tilde{t}}
	+\frac{\partial}{\partial \tilde{z}} \left(\alpha \frac{\tilde{Q}^{2}}{\tilde{A}}\right)
	+\tilde{A} \frac{\partial \tilde{p}}{\partial \tilde{z}}
	+2\frac{\alpha \tilde{Q}^{2}}{\tilde{A}} \frac{\partial \ln R_{0}(\tilde{z})}{\partial \tilde{z}}
 \textcolor{black}{+\mathcal{I}}
	=\frac{1}{\mathrm{Re}} 2 \tilde{R} \frac{\partial \tilde{u}_{z}(\tilde{R}, \tilde{z}, \tilde{t})}{\partial \tilde{r}}.
	\end{aligned}
	\label{quick_summary}
 \right.
\end{equation}
Notice that the two terms
$
\frac{\partial \tilde{p}}{\partial \tilde{z}}  \text { and } 2 \tilde{R} \frac{\partial \tilde{u}_{z}}{\partial \tilde{r}}(\tilde{R}, \tilde{z}, \tilde{t})
$
in the above equations need to be specified in terms of $\tilde{A}$ and $\tilde{Q}$ to have a closed system.

Next, to express the system (\ref{quick_summary}) in dimensional form, we define the average axial velocity as:
$$
\begin{aligned}%
	U = \frac{2}{R^2}\int_{0}^{R}ru_z dr,
\end{aligned}
$$
together with Equ.~(\ref{axial_aver_vel_nondimensional}), a straightforward derivation yields 
$
U = U_z\tilde{U}.
$
Here, we also outline the other scales we have chosen:
$$
\left\{
\begin{aligned}
& R=R_0 \tilde{R} \Rightarrow \tilde{R}=\frac{R}{R_0}, \\
& A=R_0^2 \tilde{A} \Rightarrow \tilde{A}=\frac{A}{R_0^2}, \\
& U=U_z \tilde{U} \Rightarrow \tilde{U}=\frac{U}{U_z}, \\
& t=\frac{L}{U_z} \tilde{t} \Rightarrow \tilde{t}=\frac{U_z}{L} t, \\
& z=L \tilde{z} \Rightarrow  \tilde{z}= \frac{1}{L}z,
\end{aligned}
\right.
$$
\subsubsection{Recovery of the conservation of mass in dimensional form}
By substituting the expressions of $\tilde{A}$ and $\tilde{Q}$ into Equ.~(\ref{quick_summary}.a), we get:
\begin{equation}
\begin{aligned}
& \frac{\partial \tilde{A}}{\partial \tilde{t}}+\frac{\partial \tilde{Q}}{\partial \tilde{z}}+2 \tilde{Q} \frac{\partial \ln R_0}{\partial \tilde{z}}=0 \\
\Rightarrow ~& \frac{\partial \tilde{A}}{\partial \tilde{t}}+\frac{\partial(\tilde{A} \tilde{U})}{\partial \tilde{z}}+2(\tilde{A} \tilde{U}) \frac{\partial \ln R_0}{\partial \tilde{z}}=0\\
\Rightarrow ~&\frac{L}{U_z R_0^2} \frac{\partial A}{\partial t}+L \frac{\partial\left(\frac{A U}{R_0^2 U_z}\right)}{\partial z}+2 \frac{A UL}{R_0^2 U_z} \frac{\partial \ln R_0}{\partial z}=0
\end{aligned}
\label{recovery_mass}
\end{equation}
For the second term in the equation above, we derive:
$$
\begin{aligned}
L \frac{\partial\left(\frac{A U}{R_0^2 U_z}\right)}{\partial z}&=\frac{L}{U_z} \frac{\partial\left(\frac{Q}{R_0^2}\right)}{\partial z}=\frac{L Q}{U_z} \frac{\partial \frac{1}{R_0^2}}{\partial z}+\frac{L}{R_0^2 U_z} \frac{\partial Q}{\partial z} \\
& =-\frac{2 Q L}{U_z R_0^3} \frac{\partial R_0}{\partial z}+\frac{L}{R_0^2 U_z} \frac{\partial Q}{\partial z} \\
& =-\frac{2 Q L}{U_z R_0^2} \frac{\partial \ln R_0}{\partial z} +\frac{L}{R_0^2 U_z} \frac{\partial Q}{\partial z} 
\end{aligned}
$$
Noticing that the second term on the right-hand side of the equation above cancels with the last term on the left-hand side in Equation (\ref{recovery_mass}), we derive:
$$
\frac{\partial A}{\partial t}
+ \frac{\partial Q}{\partial z} = 0.
$$

\subsubsection{Recovery of the axial momentum equation in dimensional form}
For the axial momentum equation, we derive the dimensional equation as follows:
\begin{equation}
\begin{aligned}
&\frac{\partial \widetilde{Q}}{\partial \tilde{t}}+\frac{\partial}{\partial \tilde{z}}\left(\alpha \frac{\widetilde{Q}^2}{\tilde{A}}\right)+\tilde{A} \frac{\partial \widetilde{p}}{\partial \tilde{z}}+\frac{2 \alpha \tilde{Q}^2}{\tilde{A}} \frac{\partial \ln R_0}{\partial \tilde{z}}
\textcolor{black}{ +\mathcal{I}}
=\frac{1}{R_e} 2 \tilde{R} \frac{\partial \tilde{u}_z(\tilde{R}, \tilde{z}, \tilde{t})}{\partial \tilde{r}}\\
\Rightarrow ~& \frac{L}{U_z^2 R_0^2} \frac{\partial Q}{\partial t}+\frac{L}{U_z^2} \frac{\partial}{\partial z}\left(\alpha \frac{Q^2}{A R_0^2}\right)+\frac{L A}{p_f U_z^2 R_0^2} \frac{\partial p}{\partial z}  +\frac{2 \alpha Q^2 L}{A R_0^2 U_z^2} \frac{\partial \ln R_0}{\partial z}  
\textcolor{black}{ +\mathcal{I}}
=\frac{2 v L R}{U_z^2 R_0^2} \frac{\partial u_z(R, z, t)}{\partial r},
\end{aligned}
\label{Qequation_dimensional}
\end{equation}
where $\mathcal{I}$ is given as follows:
$$
\begin{aligned}
\mathcal{I}=&-\frac{2}{35L^2}\left(\frac{\partial R_0(\tilde{z})}{\partial \tilde{z}}\right)^2
\frac{1}{\tilde{R}}
\left[
\frac{\partial \tilde{Q}^2}{\partial \tilde{z}}
-\tilde{U}^2\frac{\partial \tilde{A}^2}{\partial \tilde{z}}
+2\tilde{Q}^2\frac{\partial ln (R_0(\tilde{z})\tilde{R})}{\partial \tilde{z}}
\right]\\
=&-\frac{2}{35}\left(\frac{\partial R_0}{\partial z}\right)^2
    \frac{R_0}{R}
\left[\left(
\frac{L}{R_0^4 U_z^2}\frac{\partial Q^2}{\partial z}
-\frac{4 L Q^2}{R_0^4 U_z^2}\frac{\partial ln R_0}{\partial z}\right)
-\left(\frac{U^2 L}{U_z^2 R_0^4}\frac{\partial A^2}{\partial z}
-\frac{4 U^2 L A^2}{U_z^2 R_0^4}\frac{\partial ln R_0}{\partial z}
\right)\right.\\
&\left.+2\frac{Q^2 L}{R_0^4 U_z^2}\frac{\partial ln R}{\partial z}
\right]\\
=&-\frac{2}{35}\left(\frac{\partial R_0}{\partial z}\right)^2
    \frac{R_0}{R}
\left[
\frac{L}{R_0^4 U_z^2}\frac{\partial Q^2}{\partial z}
-\frac{U^2 L}{U_z^2 R_0^4}\frac{\partial A^2}{\partial z}
+2\frac{Q^2 L}{R_0^4 U_z^2}\frac{\partial ln R}{\partial z}
\right]\\
=&-\frac{2}{35}\left(\frac{\partial R_0}{\partial z}\right)^2
    \frac{1}{R R_0}
 \left(\frac{L}{U_z^2 R_0^2} \right)
\left[
\frac{\partial Q^2}{\partial z}
-U^2 \frac{\partial A^2}{\partial z}
+2Q^2\frac{\partial ln R}{\partial z}
\right]\\
\end{aligned}
$$
Note that the $3rd$ term 
$2Q^2\frac{\partial ln R}{\partial z}$
inside the brackets in the equation above can be derived as follows:
$$
2Q^2\frac{\partial ln R}{\partial z}
=\frac{2Q^2}{R}\frac{\partial R}{\partial z}
=\frac{Q^2}{ R^2}\frac{\partial R^2}{\partial z}
=\frac{Q^2}{ A}\frac{\partial A}{\partial z}
=\frac{Q^2}{2 A^2}\frac{\partial A^2}{\partial z},
$$
which can be combined with the $2nd$ term, allowing us to simplify the expression for $\mathcal{I}$ as follows:
\textcolor{black}{
$$
\begin{aligned}
\mathcal{I}=&-\frac{2}{35}\left(\frac{\partial R_0}{\partial z}\right)^2
    \frac{1}{R R_0}
 \left(\frac{L}{U_z^2 R_0^2} \right)
\left[
\frac{\partial Q^2}{\partial z}
-\frac{Q^2}{2A^2} \frac{\partial A^2}{\partial z}
\right]\\
=&-\frac{2}{35}\left(\frac{\partial R_0}{\partial z}\right)^2
    \frac{1}{R_0}
 \left(\frac{L}{U_z^2 R_0^2} \right)
\left[
\frac{1}{\sqrt{A}}\frac{\partial Q^2}{\partial z}
-\frac{Q^2}{A^{3/2}} \frac{\partial A}{\partial z}
\right]\\
=&-\frac{2}{35}\left(\frac{\partial R_0}{\partial z}\right)^2
    \frac{1}{R_0}
\left(\frac{L}{U_z^2 R_0^2} \right)
 \left[\sqrt{A}\frac{\partial \frac{Q^2}{A}}{\partial z}
\right]\\
=&-\frac{2}{35}\left(\frac{\partial R_0}{\partial z}\right)^2
    \tilde{R}
\left(\frac{L}{U_z^2 R_0^2} \right)
\left[ \frac{\partial}{\partial z}\left( \frac{Q^{2}}{A }\right)
\right]\\
\approx&-\frac{2}{35}\left(\frac{\partial R_0}{\partial z}\right)^2
\left(\frac{L}{U_z^2 R_0^2} \right)
\left[ \frac{\partial}{\partial z}\left( \frac{Q^{2}}{A }\right)
\right],\\
\end{aligned}
$$
given that $R=\sqrt{A}$, $\frac{R}{R_0}=\tilde{R}$, and $\tilde{R}\approx1$.
}

For the $2{nd}$ term in the Equ.~(\ref{Qequation_dimensional}) above, we have:
$$
\frac{L}{U_z^2} \frac{\partial}{\partial z}\left(\alpha \frac{Q^2}{A R_0^2}\right) 
= \frac{L}{U_z^2 R_0^2} \frac{\partial}{\partial z}\left(\alpha \frac{Q^2}{A}\right)-\frac{2 \alpha L Q^2}{A U_z^2 R_0^3} \frac{\partial R_0}{\partial z},
$$
Since $R_0$ is a function of $z$, the $2{nd}$ term above will cancel out the $4{th}$ term in Equ.~(\ref{Qequation_dimensional}). 

\subsubsection{(A, Q) system of Eqs.~(\ref{quick_summary}) in the dimensional form}\label{rectifiedAQ}
Ultimately, we obtain the (A, Q) system of Eqs.~(\ref{quick_summary}) in the dimensional form:
\begin{equation}
\left\{
\begin{aligned}
		&\frac{\partial A}{\partial t}
		+ \frac{\partial Q}{\partial z} = 0,\\
		&\frac{\partial Q}{\partial t}
		+ \frac{\partial}{\partial z}\left(\alpha \frac{Q^{2}}{A }\right)
		+\frac{A}{\rho_{f}} \frac{\partial p}{\partial z}
\textcolor{black}{-\frac{2}{35}\left(\frac{\partial R_0}{\partial z}\right)^2
\left[ \frac{\partial}{\partial z}\left( \frac{Q^{2}}{A }\right)
\right]}
		=\frac{2}{Re} R\frac{\d u_z(R, z, t)}{\d r}.
	\end{aligned}
 \right.
 \label{AQdimensional}
\end{equation}
We now introduce a correction value $\alpha_c$ to the Coriolis coefficient $\alpha$ as follows: 
\begin{equation}
    \alpha_c = -\frac{2}{35}\left(\frac{\partial R_0}{\partial z}\right)^2,
\end{equation}
noting that this coefficient arises from the stenosis of the vessel and vanishes in the case of a straight vessel. By rewriting the correction term in terms of partial derivative with respect to $z$, we obtain the following rectified (A, Q) system:
\begin{equation}
\left\{
	\begin{aligned}
		\frac{\partial A}{\partial t}
		+ \frac{\partial Q}{\partial z} &= 0,\\
		\frac{\partial Q}{\partial t}
		+ \frac{\partial}{\partial z}\left[(\alpha\textcolor{black}{+\alpha_c}) \frac{Q^{2}}{A }\right]
		+\frac{A}{\rho_{f}} \frac{\partial p}{\partial z}
		&=\frac{2}{Re} R\frac{\d u_z(R, z, t)}{\d r}\textcolor{black}{+\frac{Q^2}{A}\frac{\partial\alpha_c}{\partial z}}.
	\end{aligned}
 \right.
\end{equation}
\begin{remark}
This approximation allows one to easily reuse their legacy code for solving the proposed (A, Q) system. Another way for approximating the integration $\mathcal{I}$ is provided in \ref{section:AppendixB}. While this alternative way is more accurate in derivation, it is less concise in its final format because $\mathcal{I}$ cannot be rewritten as a correction term of any existing terms in the classical (A, Q) system. Additionally, numerical results do not show a significant advantage for this approximation. Therefore, we choose to use the first approximation approach, proposed in Sec.~\ref{rectifiedAQ}.
\end{remark}

\subsection{The structure equation to close the system}
Since the aforementioned dimensional (A, Q) system is not closed, we employ a structural model along with suitable velocity profiles to determine the pressure term $\frac{\partial p}{\partial z}$ and the velocity term $\frac{\partial u_z(R, z, t)}{\partial r}$.
We assume that the arterial walls are uniform and isotropic. A straightforward structural model option is the Koiter shell model, which focuses solely on radial displacement and disregards terms with higher derivatives. Recall the dynamic coupling condition:
\begin{equation}
\rho_{K} h \frac{\partial^{2} \eta_{r}}{\partial t^{2}}+\frac{h E}{R_{0}\left(1-\sigma^{2}\right)} \frac{\eta_{r}}{R_{0}}=-\left.\mathcal{J} \sigma \mathbf{n}\right|_{\Gamma(t)} \cdot \mathbf{e}_{r}, 
\label{membrane_equation}
\end{equation}
which indicates that the force exerted by the fluid on the membrane converts into two types of energy: kinetic energy (the first term on the left-hand side) and elastic energy (the second term on the left-hand side). The elastic energy is linked to either the stretching ($\eta_{r}>0$) or the recoil ($\eta_{r}<0$) of the membrane, characterized by its stiffness $E$, thickness $h$, and Poisson ratio $\sigma$. Meanwhile, the kinetic energy is tied to the membrane's velocity and acceleration, with a mass of $\rho_{K} h$. In the following, we transform this into non-dimensional form and make the leading-order approximation.
	
The Cauchy stress tensor $\sigma$ comprises a pressure component and a deviatoric stress term involving the symmetrized gradient of velocity. In the case of axially symmetric flow, we have that:
$$
\mathbf{D}(\mathbf{u})
=\frac{1}{2}\left[\begin{array}{ccc}
2 \frac{\partial u_{r}}{\partial r} & 0 & \frac{\partial u_{z}}{\partial r}+\frac{\partial u_{r}}{\partial z} \\
0 & 2\frac{\partial u_{r}}{\partial r}  & 0\\
\frac{\partial u_{z}}{\partial r}+\frac{\partial u_{r}}{\partial z} & 0 & 2 \frac{\partial u_{z}}{\partial z}
\end{array}\right],
$$
and the unit normal $\mathbf{n}$ to the artery wall $R(z, t)$ is given by
$$
\mathbf{n}=\frac{\left(n_{r}, n_{\theta}, n_{z}\right)}{\sqrt{n_{r}^{2}+n_{\theta}^{2}+n_{z}^{2}}}=\frac{\left(1,0,-\partial_{z} R\right)}{\sqrt{1+\left(\partial_{z} R\right)^{2}}}
=\frac{1}{\sqrt{1+\left(\partial_{z} R\right)^{2}}}\left(\mathbf{e}_{r}-\partial_{z} R \mathbf{e}_{z}\right).
$$
For the term $\left.\mathcal{J} \sigma \mathbf{n}\right|_{\Gamma(t)} \cdot \mathbf{e}_{r}$ in Equ.~(\ref{membrane_equation}), we note that the scaling factor in $\mathbf{n}$ will cancel out with the Jacobian $\mathcal{J}$, and thus we have:
	$$
	\mathcal{J}(\sigma \mathbf{n}) \cdot \mathbf{e}_{r}
 	=\mathcal{J}(-p \mathbf{I}+2 \mu \mathbf{D}(\mathbf{u})) \mathbf{n} \cdot \mathbf{e}_{r}
 	=-p +2 \mu \frac{\partial u_{r}}{\partial r}
    -\mu\left(\frac{\partial u_{z}}{\partial r}
	+\frac{\partial u_{r}}{\partial z}\right) \frac{\partial R}{\partial z}.
	$$
Recalling that $U_{r}$ and $R_{0}$ are variable with $\tilde{z}$, yields:
$$
\frac{\partial R}{\partial z}
=\frac{\tilde{R}}{L}\frac{\partial R_{0}(\tilde{z})}{\partial\tilde{z}}
+\frac{R_{0}(\tilde{z})}{L}\frac{\partial \tilde{R}}{\partial \tilde{z}}
\ \ \text{and} \ \
\frac{\partial u_{r}}{\partial z}
=-\frac{U_r(\tilde{z})\tilde{r}}{L}\frac{\partial ln R_0(\tilde{z})}{\partial\tilde{z}}\frac{\partial\tilde{u_r}}{\partial\tilde{r}}
+\frac{\tilde{u_r}}{L} \frac{\partial U_{r}(\tilde{z})}{\partial\tilde{z}}
+\frac{U_{r}(\tilde{z})}{L} \frac{\partial \tilde{u}_{r}}{\partial \tilde{z}}.
$$ 
Therefore, we derive the expression of $\mathcal{J}(\boldsymbol{\sigma} \mathbf{n}) \cdot \mathbf{e}_{r}$ in non-dimensional variables
and replace $\mathrm{Re}=\frac{\rho_f U_{z} R_{0}(z)^{2}}{\mu_f L}$ to obtain:
$$
 \begin{aligned}
{\mathcal{J}}({\boldsymbol{\sigma}} {\mathbf{n}}) \cdot \mathbf{e}_{r} 
 =&-p 
 +2\mu\frac{\partial u_{r}}{\partial r}
 -\mu\left(\frac{\partial u_{z}}{\partial r}
 +\frac{\partial u_{r}}{\partial z}\right) \frac{\partial R}{\partial z} \\
 =& -\rho_f U_z^2\tilde{p}
 + 2\mu\frac{U_r(\tilde{z})}{R_0(\tilde{z})}\frac{\partial\tilde{u}_r}{\partial\tilde{r}}
 - \mu\frac{U_z}{R_0(\tilde{z})}\frac{\partial\tilde{u}_z}{\partial\tilde{r}}\left(
	 \frac{\tilde{R}}{L}\frac{\partial R_{0}(\tilde{z})}{\partial\tilde{z}}
	 +\frac{R_{0}(\tilde{z})}{L}\frac{\partial \tilde{R}}{\partial \tilde{z}}
 \right) \\
& - \mu\left(
	 -\frac{U_r(\tilde{z})\tilde{r}}{L}\frac{\partial ln R_0(\tilde{z})}{\partial\tilde{z}}\frac{\partial\tilde{u}_r}{\partial\tilde{r}}
	 +\frac{U_r(\tilde{z})}{L}\frac{\partial\tilde{u}_r}{\partial\tilde{z}} 
	 + \frac{\tilde{u}_r}{L}\frac{\partial U_r(\tilde{z})}{\partial\tilde{z}} 
 \right)\left(
	 \frac{\tilde{R}}{L}\frac{\partial R_{0}(\tilde{z})}{\partial\tilde{z}}
	 +\frac{R_{0}(\tilde{z})}{L}\frac{\partial \tilde{R}}{\partial \tilde{z}}
 \right) \\
 =& -\rho_fU_z^2\tilde{p}
 + \rho_fU_z^2\frac{2U_r(\tilde{z})R_0(\tilde{z})}{ReU_zL}\frac{\partial\tilde{u}_r}{\partial\tilde{r}}
 - \rho_fU_z^2\frac{R_0(\tilde{z})\tilde{R}}{R_eL^2}\frac{\partial\tilde{u}_z}{\partial\tilde{r}}\frac{\partial R_0(\tilde{z})}{\partial\tilde{z}}
 - \rho_fU_z^2\frac{R_0^2(\tilde{z})}{R_eL^2}\frac{\partial\tilde{u}_z}{\partial\tilde{r}}\frac{\partial\tilde{R}}{\partial\tilde{z}}\\
& + \rho_fU_z^2\frac{U_r(\tilde{z})R_0^2(\tilde{z})\tilde{r}\tilde{R}}{R_eL^3U_z}\frac{\partial ln R_0(\tilde{z})}{\partial\tilde{z}}\frac{\partial\tilde{u}_r}{\partial\tilde{r}}\frac{\partial R_0(\tilde{z})}{\partial\tilde{z}}
 + \rho_fU_z^2\frac{U_r(\tilde{z})R_0^3(\tilde{z})\tilde{r}}{R_eL^3U_z}\frac{\partial ln R_0(\tilde{z})}{\partial\tilde{z}}\frac{\partial\tilde{u}_r}{\partial\tilde{r}}\frac{\partial\tilde{R}}{\partial\tilde{z}} \\
& - \rho_fU_z^2\frac{U_r(\tilde{z})R_0^2(\tilde{z})\tilde{R}}{R_eL^3U_z}\frac{\partial\tilde{u}_r}{\partial\tilde{z}}\frac{\partial R_0(\tilde{z})}{\partial\tilde{z}}
 - \rho_fU_z^2\frac{U_r(\tilde{z})R_0^3(\tilde{z})}{R_eL^3U_z}\frac{\partial\tilde{u}_r}{\partial\tilde{z}}\frac{\partial\tilde{R}}{\partial\tilde{z}}\\
& - \rho_fU_z^2\frac{R_0^2(\tilde{z})\tilde{u}_r\tilde{R}}{R_eL^3U_z}\frac{\partial U_r(\tilde{z})}{\partial\tilde{z}}\frac{\partial R_0(\tilde{z})}{\partial\tilde{z}}
 - \rho_fU_z^2\frac{R_0^3(\tilde{z})\tilde{u}_z}{R_eL^3U_z}\frac{\partial U_r(\tilde{z})}{\partial\tilde{z}}\frac{\partial\tilde{R}}{\partial\tilde{z}}  
 \end{aligned}
$$
By taking the leading order approximation, we have:
$$
{\mathcal{J}}({\boldsymbol{\sigma}} {\mathbf{n}}) \cdot \mathbf{e}_{r} 
\approx \rho_{f} U_{z}^{2}\left(
-\tilde{p}
- \frac{R_0(\tilde{z})\tilde{R}}{\operatorname{Re} L^2} \frac{\partial \tilde{u}_{z}}{\partial \tilde{r}} \frac{\partial R_{0}(\tilde{z})}{\partial\tilde{z}}
\right)
=  -p
-\frac{\rho_{f} U_{z} R}{2 \operatorname{Re} L} \frac{\partial u_{z}}{\partial r} \frac{\partial R_{0}^{2}(z)}{\partial z}.
$$
Assuming that
$$
\frac{\eta_{r}}{R_{0}}=\frac{R_{0}+\eta_{r}-R_{0}}{R_{0}}=\frac{R-R_{0}}{R_{0}},
$$
We can rewrite Equ.~(\ref{membrane_equation}) as:
$$
\rho_{K} h \frac{\partial^{2} R}{\partial t^{2}}+\frac{h E}{R_{0}^2\left(1-\sigma^{2}\right)}\left({R}-{R_{0}}\right)= -p
-\frac{\rho_{f} U_{z} R}{2 \operatorname{Re} L} \frac{\partial u_{z}}{\partial r} \frac{\partial R_{0}^{2}(z)}{\partial z}.
$$
\if(0)
Taking the non-dimensional form of the above equation, we have:
$$
\rho_{K} h \frac{R_{0}}{T^{2}} \frac{\partial^{2} \tilde{R}}{\partial \tilde{t}^{2}}+\frac{h E}{R_{0}\left(1-\sigma^{2}\right)}(\tilde{R}-1)=-\left.{\mathcal{J}} \tilde{\sigma} {\mathbf{n}}\right|_{\Gamma(t)} \cdot \mathbf{e}_{r},
$$
which we simplify to:
$$
\rho_{K} h \frac{R_{0}}{L} \frac{U_{z} U_{r}}{R_{0}} \frac{\partial^{2} \tilde{R}}{\partial \tilde{t}^{2}}
+\frac{h E}{R_{0}\left(1-\sigma^{2}\right)}(\tilde{R}-1)
= \rho_{f} U_{z}^{2} \tilde{p}
+ \frac{\rho_{f} U_{z} R}{2 \operatorname{Re} L} \frac{\partial u_{z}({R}, {z}, {t})}{\partial r} \frac{\partial R_{0}^{2}}{\partial z}.
$$
\fi
Ignore the term contains time derivatives, we obtain:
\if(0)
$$
\frac{h E}{R_{0}\left(1-\sigma^{2}\right)}(\tilde{R}-1)
= p
+ \frac{\rho_{f} U_{z} R}{2 \operatorname{Re} L} \frac{\partial u_{z}}{\partial r} \frac{\partial R_{0}^{2}}{\partial z},
$$
or
\fi
$$
\frac{h E}{R_{0}^2\left(1-\sigma^{2}\right)}\left({R}-R_{0}\right) 
-\rho_{f} \nu R \frac{\partial u_{z}}{\partial r} \frac{\partial \ln R_{0}}{\partial z}=p-p_{\mathrm{ext}},
$$
where $\nu=\frac{\mu_f}{\rho_f}$. The dynamic coupling condition states that the elastodynamics of the artery wall are driven by the change in pressure. What is left is to define the viscous term $\frac{\partial u_z}{\partial r}$ in terms of ${A}$ and ${Q}$, and to determine $\alpha$.  For the rest of the manuscript, we revert to the dimensional representation of our (A, Q) system, as no further asymptotic analysis will be pursued.
	
\subsection{The viscous term and Coriolis factor $\alpha$}
To be able to capture the velocity profiles for a range of Womersley numbers, a more general ad-hoc closure assumption on the velocity profile can be used, which takes the form:
$$
u_{z}=\frac{\gamma+2}{\gamma} U
\left[1-\left(\frac{r}{R}\right)^{\gamma}\right].
$$
Namely, we are assuming polynomial behavior of degree $\gamma .$ The larger the $\gamma,$ the flatter the velocity profile.
With this $\gamma,$ one can show using direct calculations, that the Coriolis coefficient $\alpha$ and the viscous term in the momentum equation become:
$$
\gamma=\frac{2-\alpha}{\alpha-1}
$$
and
$$
\nu\left[2 R \frac{\partial u_{z}}{\partial r}(R, z, t)\right]=-2(\gamma+2) \nu U=-2(\gamma+2) \nu \frac{Q}{A}.
$$
For all of our simulations, we chose $\gamma = 9$ and $\alpha = 1.1$.
	
\section{The closed $1D$ reduced (A, Q) system in dimensional form}	
The complete extended $1D$ system in dimensional form is as follows:
\begin{equation}
\left\{
\begin{aligned}
&\frac{\partial A}{\partial t}
+ \frac{\partial Q}{\partial z} = 0,\\
&\frac{\partial Q}{\partial t}
+ \frac{\partial}{\partial z}\left[(\alpha\textcolor{black}{+\alpha_c}) \frac{Q^{2}}{A }\right]
+\frac{A}{\rho_{f}} \frac{\partial p}{\partial z}
=-2(\gamma+2) \nu \frac{Q}{A}
+\frac{Q^2}{A}\frac{\partial\alpha_c}{\partial z},\\
&p=
p_{\text {ext }}
+ \frac{h E}{R_{0}^2\left(1-\sigma^{2}\right)}\left({R}-{R_{0}}\right)
+ (\gamma+2) \rho_{f} \nu \frac{Q}{ A } \frac{\partial \ln R_{0}}{\partial z}.
\end{aligned}
\label{AQsystem2}
\right.
\end{equation}
where $\gamma=-\frac{\alpha-2}{\alpha-1}$.

To express the system (\ref{AQsystem2}) in conservative form, we separate Equ.~(\ref{AQsystem2}c) into two parts:
\begin{equation}
p_{1}(A)
=p_{\mathrm{ext}}+\frac{h E}{R_{0}^2\left(1-\sigma^{2}\right)}\left({R}-{R_{0}}\right),
	\label{p1}
	\end{equation}
	and
\begin{equation}\quad 
p_{2}(A, Q)
=(\gamma+2) \rho_{f} \nu \frac{Q}{ A } \frac{\partial \ln R_{0}}{\partial z}.
\label{p2}
\end{equation}
Then, we reconstruct the term $\frac{A}{\rho_{f}} \frac{\partial p_{1}}{\partial z}$
in pseudo-conservative form:
	\begin{equation}
	\begin{aligned}
	&\frac{A}{\rho_{f}} \frac{\partial p_{1}}{\partial z}
	=\frac{\partial}{\partial z}\left(\frac{h E R^{3}}{3 \rho_{f}\left(1-\sigma^{2}\right) R_{0}^{2}}\right)-\frac{4 h E R^{3}}{3 \rho_{f}\left(1-\sigma^{2}\right) R_{0}^{3}} \frac{\partial R_{0}}{\partial z}+\frac{h E A}{\rho_{f}\left(1-\sigma^{2}\right) R_{0}^{2}} \frac{\partial R_{0}}{\partial z}.
	\label{p1p2_conser}
	\end{aligned}
	\end{equation}
Thus, we have:
\begin{equation}
	\begin{aligned}
	\frac{\partial Q}{\partial t}+\frac{\partial}{\partial z}\left[(\alpha+\alpha_c) \frac{Q^{2}}{A}
	+\frac{h E A^{\frac{3}{2}}}{3 \rho_{f}\left(1-\sigma^{2}\right) R_{0}^{2}}\right]
	+\frac{A}{\rho_{f}} \frac{\partial p_{2}}{\partial z}
	=-2(\gamma+2) \nu \frac{Q}{A}\\
	+\frac{ h E A^{\frac{3}{2}}}{3 \rho_{f}\left(1-\sigma^{2}\right) R_{0}^{3}} \frac{\partial R_{0}}{\partial z}
	-\frac{h E A}{\rho_{f}\left(1-\sigma^{2}\right) R_{0}^{2}} \frac{\partial R_{0}}{\partial z}+\frac{Q^2}{A}\frac{\partial\alpha_c}{\partial z},
	\label{1D_Q_conser}
	\end{aligned}
\end{equation}
where $\frac{A}{\rho_f}\frac{\partial p_2}{\partial z}=(\gamma+2)\nu \left[ \left(  \frac{dQ}{dz} - \frac{Q}{A} \frac{dA}{dz} \right) \cdot \frac{\partial \ln(R_0)}{\partial z} + \frac{Q}{A} \cdot \frac{\partial^2 \ln(R_0)}{\partial z^2} \right]$ will be moved to the RHS.

\subsection{Conservative form and Riemann Invariants}
\label{section:RI}
We express Equs.~(\ref{AQsystem2}) into the conservative form:
	$$
	\frac{\partial U}{\partial t}+\frac{\partial F(U)}{\partial z}=S(U),
	$$
	where
	$$
	U=\left[\begin{array}{l}
	A \\
	Q
	\end{array}\right],
	\quad
	F(U)=\left[\begin{array}{c}
	Q \\
	\frac{(\alpha+\alpha_c) Q^{2}}{A}+\frac{h E A^{\frac{3}{2}}}{3 \rho_{f}\left(1-\sigma^{2}\right) R_{0}^{2}}
	\end{array}\right],
	$$
	and
$$
S(U)=
\left[\begin{array}{c}
0 \\
-2(\gamma+2) \nu \frac{Q}{A}
+\frac{4 h E A^{\frac{3}{2}}}{3 \rho_{f}\left(1-\sigma^{2}\right) R_{0}^{3}} \frac{\partial R_{0}}{\partial z}
-\frac{h E A}{\rho_{f}\left(1-\sigma^{2}\right) R_{0}^{2}} \frac{\partial R_{0}}{\partial z}-\frac{A}{\rho_{f}} \frac{\partial p_{2}}{\partial z}+\frac{Q^2}{A}\frac{\partial\alpha_c}{\partial z}
\end{array}\right].
$$
As for the quasilinear form
$U_{t}+F^{\prime}(U) U_{z}=S(U)$, $F^{\prime}(U)$ has the following expression:
$$
F^{\prime}(U)=\left[\begin{array}{cc}
0 & 1\\
-(\alpha+\alpha_c) \frac{Q^{2}}{A^{2}}+\frac{h E \sqrt{A}}{2 \rho_{f}\left(1-\sigma^{2}\right) R_{0}^{2}}
& 2 (\alpha+\alpha_c) \frac{Q}{A}
\end{array}\right]
$$
Denote $B=F^{\prime}(U)$ and calculate its eigenvalues:
$$
\lambda^{2}-\operatorname{Tr}(B) \lambda+\operatorname{det}(B)=0,
$$
we have:
$$
\lambda_{1,2}=\frac{\operatorname{Tr}(B) \mp \sqrt{\operatorname{Tr}(B)^{2}-4 \operatorname{det}(B)}}{2},
$$
where
$$
\operatorname{Tr}(B)=2 (\alpha+\alpha_c) \frac{Q}{A},
\quad
\operatorname{det}(B)=(\alpha+\alpha_c) \frac{Q^{2}}{A^{2}}-\frac{h E \sqrt{A}}{2 \rho_{f}\left(1-\sigma^{2}\right) R_{0}^{2}}.
$$
Therefore, we obtain the two eigenvalues:
$$
\begin{aligned}
\lambda_{1,2}=(\alpha+\alpha_c) \frac{Q}{A}
\mp \sqrt{\left[(\alpha+\alpha_c) \frac{Q}{A}\right]^{2}
	-(\alpha+\alpha_c) \frac{Q^{2}}{A^{2}}+\frac{h E \sqrt{A}}{2 \rho_{f}\left(1-\sigma^{2}\right) R_{0}^{2}}}.
	\end{aligned}
	$$
and the right eigenvectors:
$$
I_{1}=\left[\begin{array}{c}
1 \\
\lambda_{2}
\end{array}\right],
\quad
I_{2}=\left[\begin{array}{c}
1 \\
\lambda_{1}
\end{array}\right].
$$
Notice that regions near both inlet and outlet are straight, causing the correction value $\alpha_c$ to vanish in these regions. To obtain the Riemann invariants for rectified (A, Q) system, we consider a flat velocity profile closure, specifically $\alpha=1$, as this allows for the explicit derivation of the Riemann invariants $w_i ~(i = 1, 2)$. In this case, we have:
$$
I_{1}=\left[\begin{array}{c}
	1\\
	\frac{Q}{A} - \sqrt{\frac{hE\sqrt A}{2\rho_f(1-\sigma^2)R^2_0}} 
	\end{array}\right],
	\quad
I_{2}=\left[\begin{array}{c}
	1 \\
	\frac{Q}{A} + \sqrt{\frac{hE\sqrt A}{2\rho_f(1-\sigma^2)R^2_0}}
	\end{array}\right].
$$
Let $\beta = \frac{hE}{2\rho_f(1-\sigma^2)R_0^2}$. In characteristic form the system is
$$
	\frac{dA}{1}=\frac{dQ}{\frac{Q}{A}\mp\beta^{\frac{1}{2}}A^{\frac{1}{4}}}=\text{const}.
$$
By multiplying an integrating factor $A$, we obtain:
$$
\begin{aligned}
	A\left[\left(\frac{Q}{A}\mp\beta^{\frac{1}{2}}A^{-\frac{3}{4}}\right)dA-\frac{1}{A}dQ\right]=0\\
	\Rightarrow A\left[d\left(-\frac{Q}{A}\right)\mp\beta^{\frac{1}{2}}A^{-\frac{3}{4}}dA\right]=0\\
	\Rightarrow Ad\left(-\frac{Q}{A}\mp4\beta^{\frac{1}{2}}A^{\frac{1}{4}}\right)=0,
\end{aligned}
$$
which leads to:
\begin{equation}
	w_{1, 2} = -\frac{Q}{A}\mp4\sqrt{\frac{hE\sqrt{A}}{2\rho_f(1-\sigma^2)R_0^2}},
\end{equation}
which are the Riemann invariants.	
	
\subsection{A simple case considering isotropic-homogeneous elastic structure}
In a simple case, we assume the structure is isotropic and homogeneous elastic. We can rewrite the following equation in terms of displacement:
	\begin{equation}p_{1}(A)=p_{\mathrm{ext}}+\frac{h E}{R_{0}^2\left(1-\sigma^{2}\right)}\left({R}-{R_{0}}\right)=p_{\mathrm{ext}}+C_0 \eta,
	\end{equation}
where $C_0=\frac{h E}{R_{0}^2\left(1-\sigma^{2}\right)}$. In our model of blood flow through an artery with stenosis, $R_0$ denotes the reference radius of the artery, which varies with $z$. However, in our numerical simulation, we will set $R_0$ in the expression of $C_0$ as a constant, denoting it as $R_0^*$. $R_0^*$ is different from the reference radius $R_0$. Consequently, $C_0$ will be constant throughout the artery and will not depend on $z$. Therefore, the Equ.~(\ref{p1p2_conser}) is changed to
	$$
	\frac{A}{\rho_{f}} \frac{\partial p_{1}}{\partial z}
	=\frac{\partial}{\partial z}\left(\frac{h E R^{3}}{3 \rho_{f}\left(1-\sigma^{2}\right) (R_{0}^{*})^2}\right)
	-\frac{h E A}{\rho_{f}\left(1-\sigma^{2}\right) (R_{0}^*)^{2}} \frac{\partial R_{0}}{\partial z}.
	$$
	And the flux term and RHS term are:
\begin{equation}
	\begin{aligned}
	&F(U)=\left[\begin{array}{c}
	Q \\
	\frac{(\alpha+\alpha_c) Q^{2}}{A}+\frac{h E A^{\frac{3}{2}}}{3 \rho_{f}\left(1-\sigma^{2}\right) (R_{0}^*)^{2}}
	\end{array}\right],\\
	&S(U)=\left[\begin{array}{c}
		0 \\
		-2(\gamma+2) \nu \frac{Q}{A}
		+\frac{h E A}{\rho_{f}\left(1-\sigma^{2}\right) (R_{0}^*)^{2}}         \frac{\partial R_{0}}{\partial z}
			-\frac{A}{\rho_f}\frac{\partial p_2}{\partial z}
   +\frac{Q^2}{A}\frac{\partial\alpha_c}{\partial z}
		\end{array}\right].
	\end{aligned}
	\label{FUSU_star}
\end{equation}

\subsection{The DG formulation to solve the (A, Q) system}
We solve the (A, Q) system in conservative form with flux and RHS defined as in Equ.~(\ref{FUSU_star}).
We divide the domain $[0, L]$ uniformly into $N$ subintervals $ \left\{I_{i}=[z_{i}, z_{i+1}]\right\}_{i=0}^{N}$.
Let $\mathbb{P}^{k}(I_{i})$ be the space of polynomials of degree $k$ on the interval $I_{i}$. The approximation space is $\mathbb{V}_{h}^{k}:=\left\{\phi:\phi \rvert_{i} \in \mathbb{P}^{k}(I_{i}) \right\}$.
We define the notation for traces of a function $\phi:[0, L] \rightarrow \mathbb{R}$ to the interior boundaries of the intervals:
$$
\left.\phi^{\pm}\right|_{z_{i}}:=\lim _{\varepsilon \rightarrow 0 \text { and } \varepsilon>0} \phi\left(z_{i} \pm \varepsilon\right), \quad \text{for}~  i=1, \ldots, N
$$
Then the semi-discretized (A, Q) problem is to find $\mathbf{U}_{h} \in \mathbb{V}_{h}^{k} \times \mathbb{V}_{h}^{k}$ satisfying:
\begin{equation}
\begin{aligned}
&\int_{I_{i}} \frac{\partial \mathbf{U}_{h}}{\partial t} \cdot \Phi_{h}=\int_{I_{i}} \mathbf{F}\left(\mathbf{U}_{h}\right) \cdot \frac{d \Phi_{h}}{d z}+\int_{I_{i}} \mathbf{S}\left(\mathbf{U}_{h}\right) \cdot \Phi_{h} \\
&-\mathbf{F}^*\left(\mathbf{U}_{h}\right)|_{z_{i+1}} \cdot {\Phi}^{-}|_{z_{i+1}}+\mathbf{F}^*\left(\mathbf{U}_{h}\right)|_{z_{i}} \cdot \Phi^{+}|_{z_{i}}, \quad \text{for}~ i=0, \ldots, N,
\label{weakformulation}
\end{aligned}
\end{equation}
for all $\Phi_{h} \in \mathbb{V}_{h}^{k} \times \mathbb{V}_{h}^{k}$.
Let $\mathbf{F}^{*}\left(\mathbf{U}_{h}\right)$ denote the numerical flux. Before we write out the flux expression, we first introduce some notation:	
$$
\begin{aligned}
	&\left.\left\{\mathbf{F}\left(\mathbf{U}_{h}\right)\right\}\right|_{z_{i}} =\frac{\mathbf{F}\left(\mathbf{U}_{h}^{+} |_{ z_{i}}\right)+\mathbf{F}\left(\mathbf{U}_{h}^{-} |_{ z_{i}}\right)}{2},\\
	&\left[[\mathbf{U}_{h}\right]]|_{z_{i}}:=\left.\mathbf{U}_{h}^{+}\right|_{z_{i}}-\left.\mathbf{U}_{h}^{-}\right|_{z_{i}}, \quad \text{for}~ i=1, \ldots, N.
\end{aligned}
$$
The local Lax-Friedrichs flux is given by:
$$	\left.\mathbf{F}^*\left(\mathbf{U}_{h}\right)\right|_{z_{i}}=\left.\left\{\mathbf{F}\left(\mathbf{U}_{h}\right)\right\}\right|_{z_{i}}-\frac{1}{2} \max_{z_{i}} (\lambda)[[\mathbf{U}_{h}]]|_{ z_{i}}, \quad \text{for}~ i=0, \ldots N
$$	
For the time discretization, we choose the $3$rd order Runge-Kutta scheme:
$$
\begin{aligned}
&\mathbf{U}_{h}^{n+\frac{1}{3}}=\mathbf{U}_{h}^{n}+\Delta t M^{-1} R H S(\mathbf{U}_{h}^{n}),\\
&\mathbf{U}_{h}^{n+\frac{2}{3}}=\frac{3}{4} \mathbf{U}_{h}^{n}+\frac{1}{4} \mathbf{U}_{h}^{n+\frac{1}{3}}+\frac{\Delta t}{4}M^{-1} R H S(\mathbf{U}_{h}^{n+\frac{1}{3}}),\\
&\mathbf{U}_{h}^{n+1}=\frac{1}{3} \mathbf{U}_{h}^{n}+\frac{2}{3} \mathbf{U}_{h}^{n+\frac{2}{3} }+\frac{2\Delta t}{3}M^{-1} R H S(\mathbf{U}_{h}^{n+\frac{2}{3} }),
\end{aligned}
$$
where $RHS$ is the R.H.S of Equ.~(\ref{weakformulation}) and $M$ denotes the mass matrix.

\subsection{Revisiting the radial momentum equation}
After we solve the (A, Q) system, we can also retrieve information about the fluid velocity in the radial direction. Recall our leading order momentum equation Equ. (\ref{Rad_Moment}), which is inconsistent as written. To fix this, we will keep higher order terms. We start from the radial momentum equation as follows by dropping the inertial and longitudinal diffusion terms, which we know a priori are order $\eps^4$, so we keep the order $\eps$, $\eps^2$, and $\eps^3$ terms:
	$$
	\rho_{f}\left(u_{r} \frac{\partial u_{r}}{\partial r}+u_{z} \frac{\partial u_{r}}{\partial z}\right)-\mu_{f}\left(\frac{\partial^{2} u_{r}}{\partial r^{2}}+\frac{1}{r} \frac{\partial u_{r}}{\partial r}-\frac{u_{r}}{r^{2}}\right)+\frac{\partial p}{\partial r}=0.
	$$
	Then, we express the above equation in non-dimensional form:
	$$
	\frac{U_{r} \tilde{u}_{r}}{R_{0}} \frac{\partial U_{r} \tilde{u}_{r}}{\partial \tilde{r}}+\frac{U_{z} \tilde{u}_{z}}{L} \frac{\partial U_{r} \tilde{u}_{r}}{\partial \tilde{z}}+\frac{1}{\rho R_{0}} \frac{\partial \rho U_{z}^{2} \tilde{p}}{\partial \tilde{r}}=\frac{\nu}{R_{0}^{2}} \frac{\partial^{2} U_{r} \tilde{u}_{r}}{\partial \tilde{r}^{2}}+\frac{\nu}{R_{0}^{2} \tilde{r}} \frac{\partial U_{r} \tilde{u}_{r}}{\partial \tilde{r}}-\frac{\nu U_{r} \tilde{u}_{r}}{R_{0}^{2} \tilde{r}^{2}}.
	$$
	Multiplying by $R_{0}/U_{z}^{2}$ and take into account that:
	$$
	\nu =\frac{U_{z} R_{0}^{2}}{\operatorname{Re} L}, \quad
	U_{r} =\frac{R_{0} U_{z}}{L}, \quad
	u_{z}=2 U_{z}\left(1-\frac{r^{2}}{R^{2}}\right).
	$$
	We obtain:
	$$
	\frac{U_{r}^{2} \widetilde{u}_{r}}{U_{z}^{2}} \frac{\partial \widetilde{u}_{r}}{\partial \widetilde{r}}+\frac{\widetilde{u}_{z} \widetilde{u}_{r} R_{0}}{L^{2}} \frac{\partial R_{0}}{\partial \tilde{z}}+\frac{\partial \widetilde{p}}{\partial \widetilde{r}}=\frac{R_{0}}{\operatorname{Re} L U_{z}} \frac{\partial^{2} U_{r} \widetilde{u}_{r}}{\partial \widetilde{r}^{2}}+\frac{R_{0}}{\operatorname{Re} L U_{z} \widetilde{r}} \frac{\partial U_{r} \widetilde{u}_{r}}{\partial \tilde{r}}-\frac{\nu U_{r} \widetilde{u}_{r}}{U_{z}^{2} R_{0} \widetilde{r}^{2}}.
	$$
Next, we multiply the above equation by $\frac{U_{z}^{2}}{R_{0}}$ to return to its dimensional form:
	$$
	u_{r} \frac{\partial u_{r}}{\partial r}+\frac{u_{r} u_{z}}{R_{0}} \frac{\partial R_{0}}{\partial z}+\frac{1}{\rho_{f}} \frac{\partial p}{\partial r}=\frac{U_{r}}{\operatorname{Re}} \frac{\partial^{2} u_{r}}{\partial r^{2}}+\frac{R_{0} U_{r}}{\operatorname{Re} r} \frac{\partial u_{r}}{\partial r}-\frac{\nu \operatorname{Re} u_{r}}{U_{r} r^{2}}.
	$$
Simplifying the above equation by multiplying $\frac{r Re}{U_r}$, we get:
	$$
	r \frac{\partial^{2} u_{r}}{\partial r^{2}}+\left(R_{0}-\frac{u_{r} r \operatorname{Re}}{U_{r}}\right) \frac{\partial u_{r}}{\partial r}+\left(\frac{4 U_{z} r}{R}-\frac{2 u_{r} r \operatorname{Re} U_{z}}{R_{0} U_{r}}\left(1-\frac{r^{2}}{R^{2}}\right)\right) \frac{\partial R_{0}}{\partial z}+\frac{R_{0} u_{r}}{r}=0,
	$$
which is an second-order ordinary differential equation in terms of $u_{r}(r)$, with boundary conditions $u_{r}(0)=0$ and $u_{r}(R)=\frac{\partial R}{\partial t}$. The equation can be solved using standard ODE solvers.

\section{Numerical results}
\subsection{Comparison of the established $1D$ reduced model, our extended $1D$ model, and the full $3D$ compliant vessel simulation under steady-state condition.}
\begin{figure}[ht]
  \centering
 \begin{overpic}[width=0.88\textwidth,grid=false,tics=10]{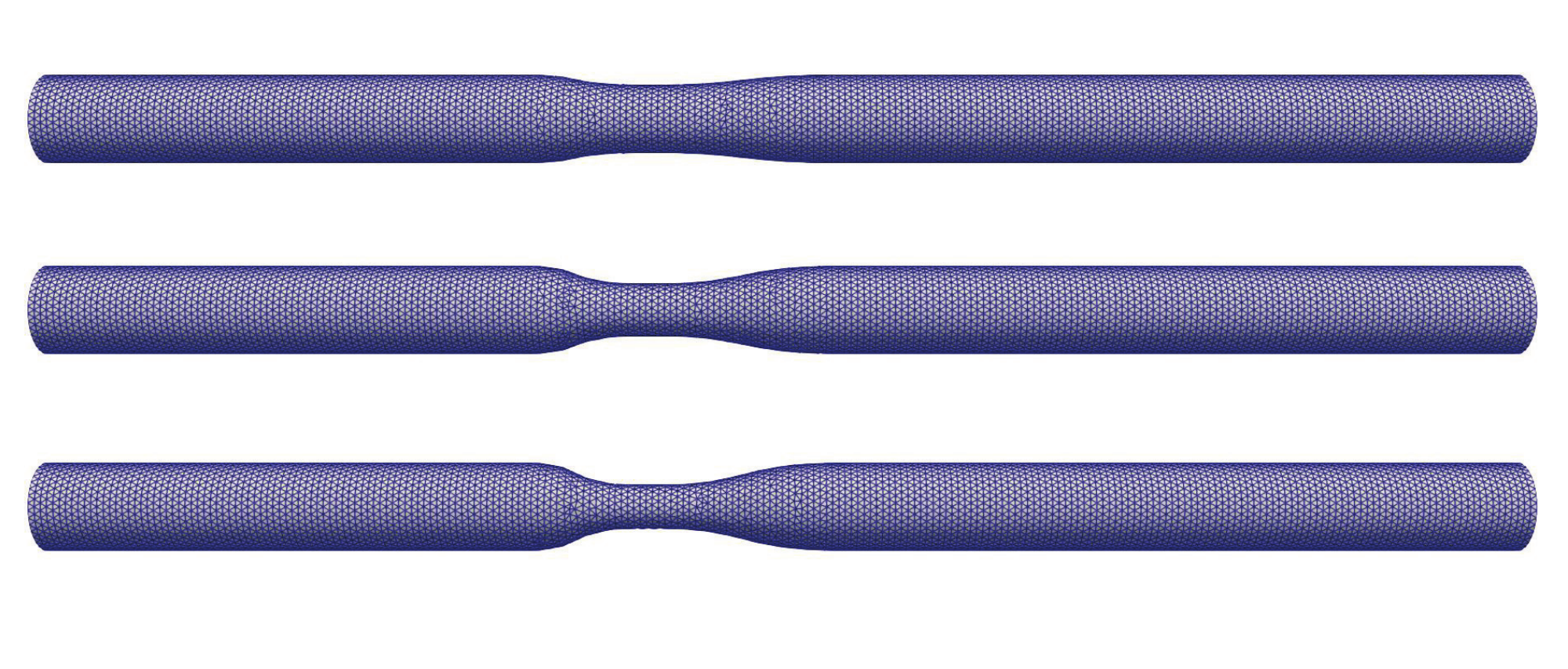}
      \put(-5,22){{Inlet}}
      \put(99,22){{Outlet}}
      \put(10,28){{23\% stenosis: ~ $R_0 =R_{max}-(R_{max}-R_{min})*e^{-50( z-3.4+0.95e^{-0.5*(z-2.5)^2 })^4}$}}
      \put(10,16){{40\% stenosis:  ~ $R_0 =R_{max}-0.4R_{max}*e^{-50( z-3.4+0.95e^{-0.5*(z-2.5)^2 })^4}$}}
       \put(10,3){{50\% stenosis:   ~ $R_0 =R_{max}-0.5R_{max}*e^{-50( z-3.4+0.95e^{-0.5*(z-2.5)^2 })^4}$}}
    \end{overpic}
  \caption{The $3D$ mesh showing the reference domain of stenotic arteries along with their corresponding geometric equations. The mesh size is around 100k tetrahedron elements.}
  \label{Fig.geo}
\end{figure}
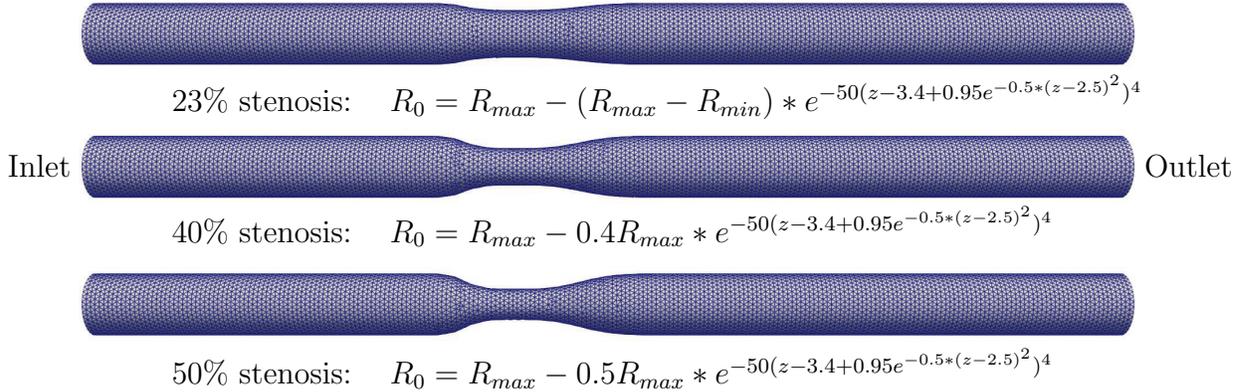

To validate our $1D$ extended model, we considered the full $3D$ fluid-structure interaction (FSI) problem of blood flow through a series of cylindrical arteries with varying radii, representing the stenosis severity from $23\%$ to $50\%$, as illustrated in Figure.~\ref{Fig.geo}.  
The artery, with a length $L$ and radius $R_0(z)$, is assumed to be an isotropic and homogeneous elastic structure and is modeled by the membrane model:
\begin{equation}
\begin{array}{ll}
\rho_{m} h \frac{\partial^{2} \eta_{r}}{\partial t^{2}}+C_{0} \eta_{r}
=f_{r} & \text { on } \Gamma^{s} \times(0, T)
\end{array}
\end{equation}
where $\eta_{r}$ denotes the radial displacement and
$$
\begin{array}{ll}
C_{0}=\frac{h}{R_0^{2}}\left(\frac{2 \mu_{m} \lambda_{m}}{\lambda_{m}+2 \mu_{m}}+2 \mu_{m}\right)\left(1+\frac{h^{2}}{12 R_0^{2}}\right),
\end{array}
$$
in which the coefficients $\mu_{m}$ and $\lambda_{m}$ are the Lam\'e coefficients, which are dependent on the Poisson ratio $\nu$ and Young's modulus $E$. The values of all parameters are provided in the tables \ref{Table1}. The artery is assumed to be clamped at both ends while allowing radial deformation. We imposed a Poiseuille flow condition of maximum velocity of $45~cm/s$ at the inlet, and a zero stress condition $\boldsymbol{\sigma} \boldsymbol{n}|_{\Gamma_{out}}= \mathbf{0}$ at the outlet.
\begin{table}[htbp]
\begin{center}\caption{Structural parameters for the test of $3D$ compliant vessel with $1D$ stent.}
\resizebox{\textwidth}{!}{
     \begin{tabular}{  l  l  l  l }
		\hline
		 Blood density (g/cm$^3$) & $\rho_f = 1.055$  & Blood Viscosity ($cm^2/s$)&  $\mu_f = 0.04 $   \\
		  Artery Length (cm)& $L = 6$  & Artery radius  (cm)&  $R_{max} ~(R_{min})= 0.18~(0.1394) $  \\
		 Artery density (g/cm$^3$)& $\rho_s = 1.055$  & Artery thickness (cm) &  $h = 0.06 $  \\
         Artery Poisson ratio  &   $\nu = 0.5$               & Artery Young modulus ($dynes/cm^2$)&   $E = 5.02 \times 10^6$   \\
         \hline
		\label{Table1}
		\end{tabular}
		}
\end{center}
\end{table}
The numerical solver developed in \cite{BUKAC2019679} has been utilized to solve the FSI problem and generate numerical solutions for comparison. A steady-state $3D$ solutions at $T=1\ s$ are chosen as a benchmark for validating the $1D$ solutions. Specifically, for $3D$ results, we computed the averaged cross-sectional velocity and pressure to serve as our validation data. 
\begin{figure}[htbp]
	\center
 \subfigure[23\% stenosis]{ \includegraphics[width=0.55\columnwidth]{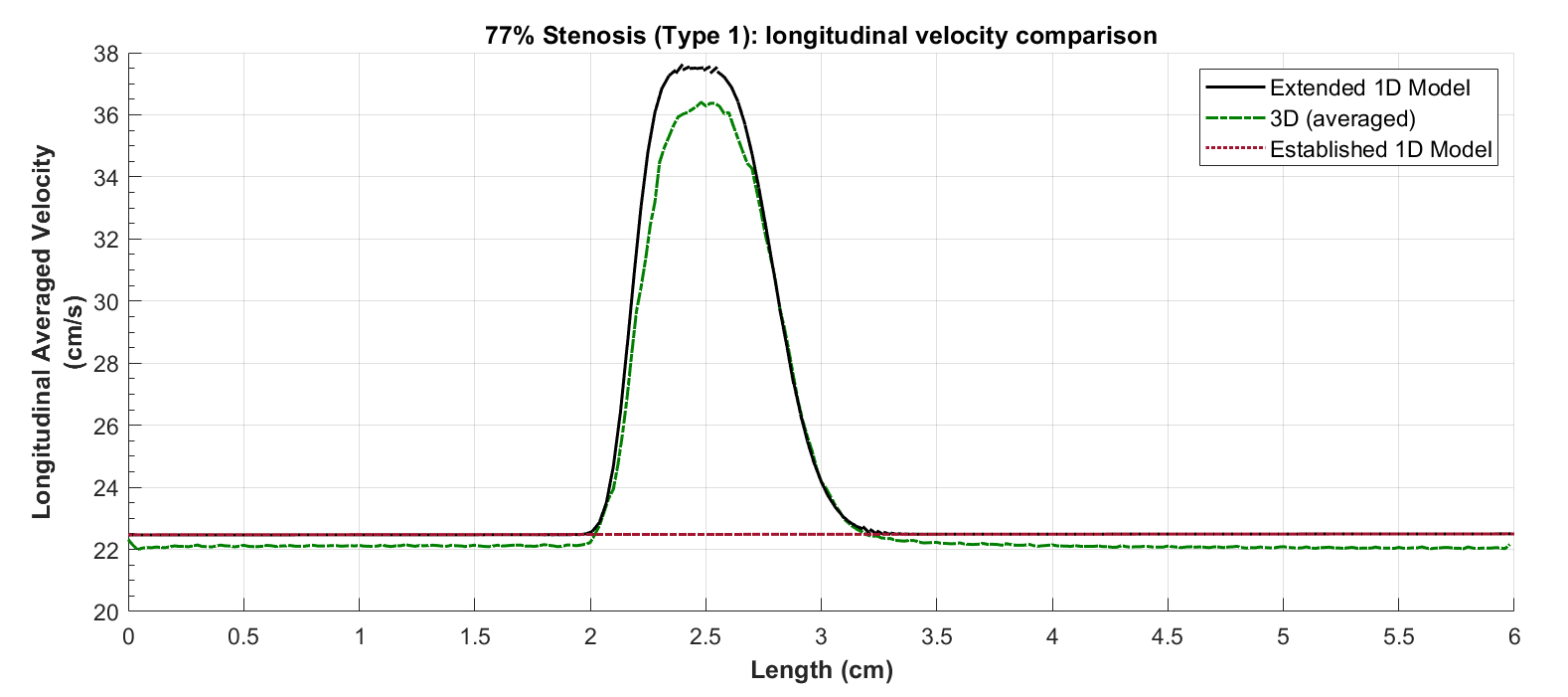}}
  \subfigure[40\% stenosis]{   \includegraphics[width=0.55\columnwidth]{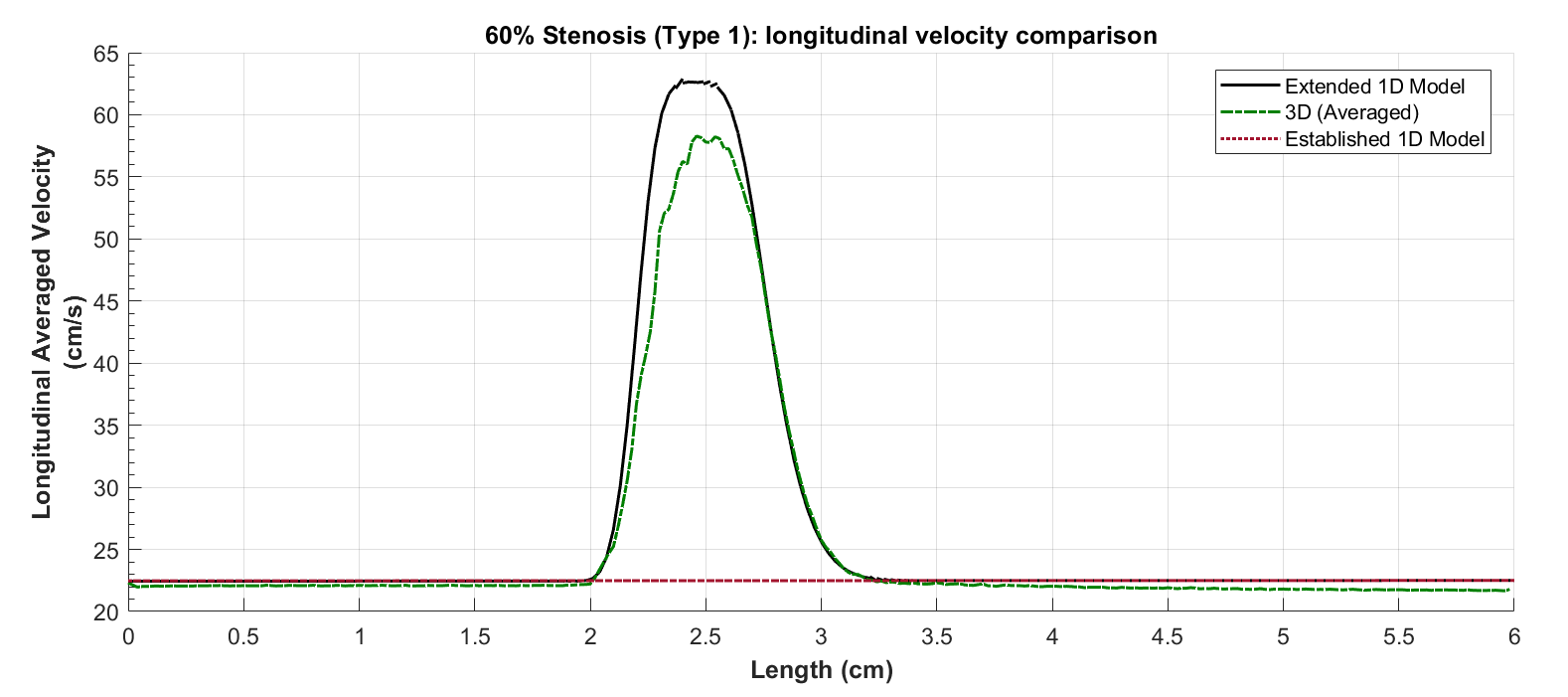}}
  \subfigure[50\% stenosis]{ \includegraphics[width=0.55\columnwidth]{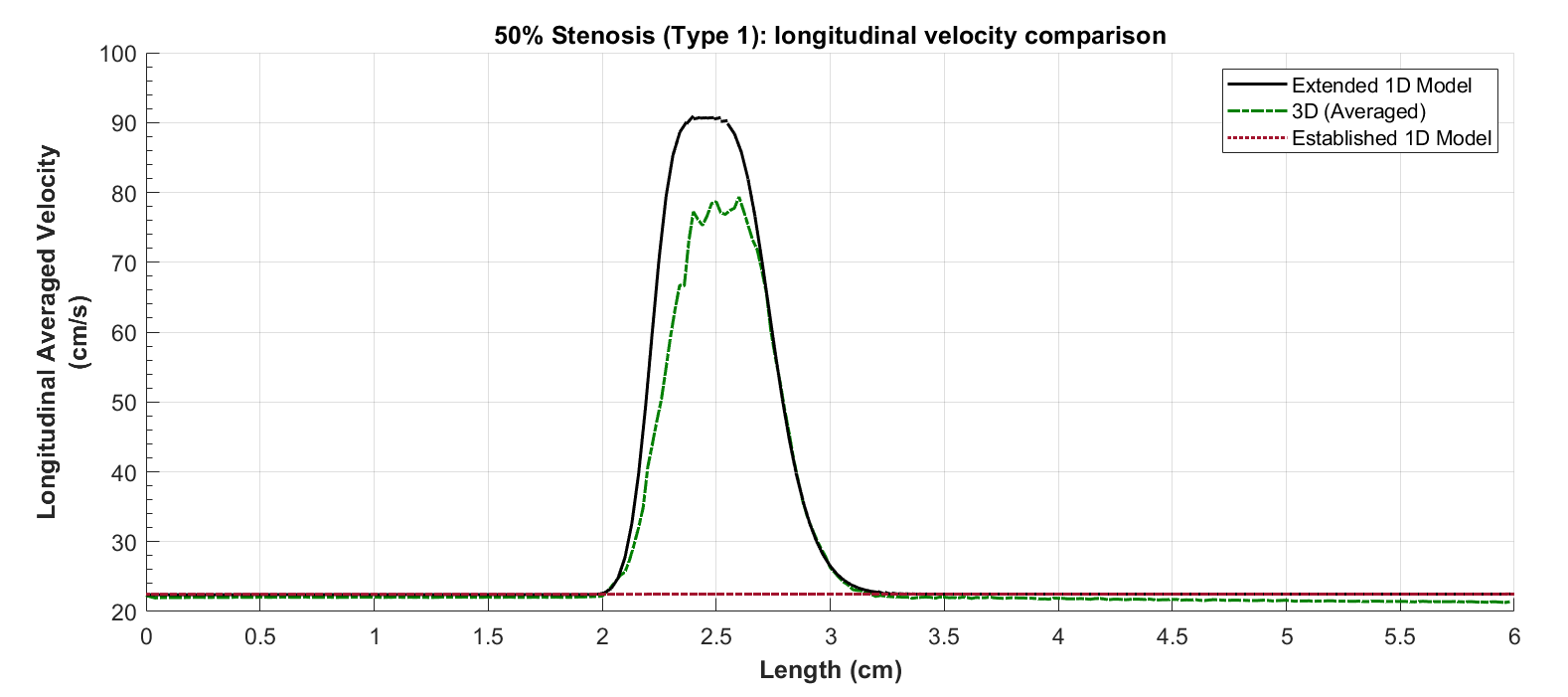}}
\caption{The longitudinal velocity $u_z$ obtained from our extended $1D$ model (in green) is compared with the established $1D$ model (in red) and the full $3D$ model (in black) results across all cases. Our extended $1D$ model shows good agreement with the $3D$ results, while the established $1D$ model fails.}
	\label{Fig.uz}
\end{figure}

Figure \ref{Fig.uz} shows the longitudinal velocity profiles, containing the results obtained from the established $1D$ model, the extended $1D$ model, and the $3D$ simulation, all superimposed for direct comparison. We observed that the extended $1D$ model successfully captures the trend of longitudinal velocity for all the cases, with a maximum error within 10\%. In contrast, the established model does not capture the correct velocity profile, as indicated by the red dotted lines. We notice some differences in the magnitude of the velocity near the stenosis between the extended $1D$ model results and the $3D$ results, which we attribute primarily to two factors. First, in the extended $1D$ model, the assumption of a parabolic velocity profile even near the entry of the stenotic region may result in an overestimate of the incoming velocity magnitude. In contrast, in the $3D$ case, we reported the averaged cross-sectional longitudinal velocity. Near the exit of the stenotic region, the outflow diverges radially due to the channel expansion, which may not be well captured by the extended $1D$ model, thus contributing to the overestimation of the extended $1D$ results. These inferences are also supported by the observation that the overall discrepancy between the extended $1D$ and $3D$ results increases with the severity of the stenosis.
\begin{figure}[htbp]
	\center
 \subfigure[23\% stenosis]{ \includegraphics[width=0.52\columnwidth]{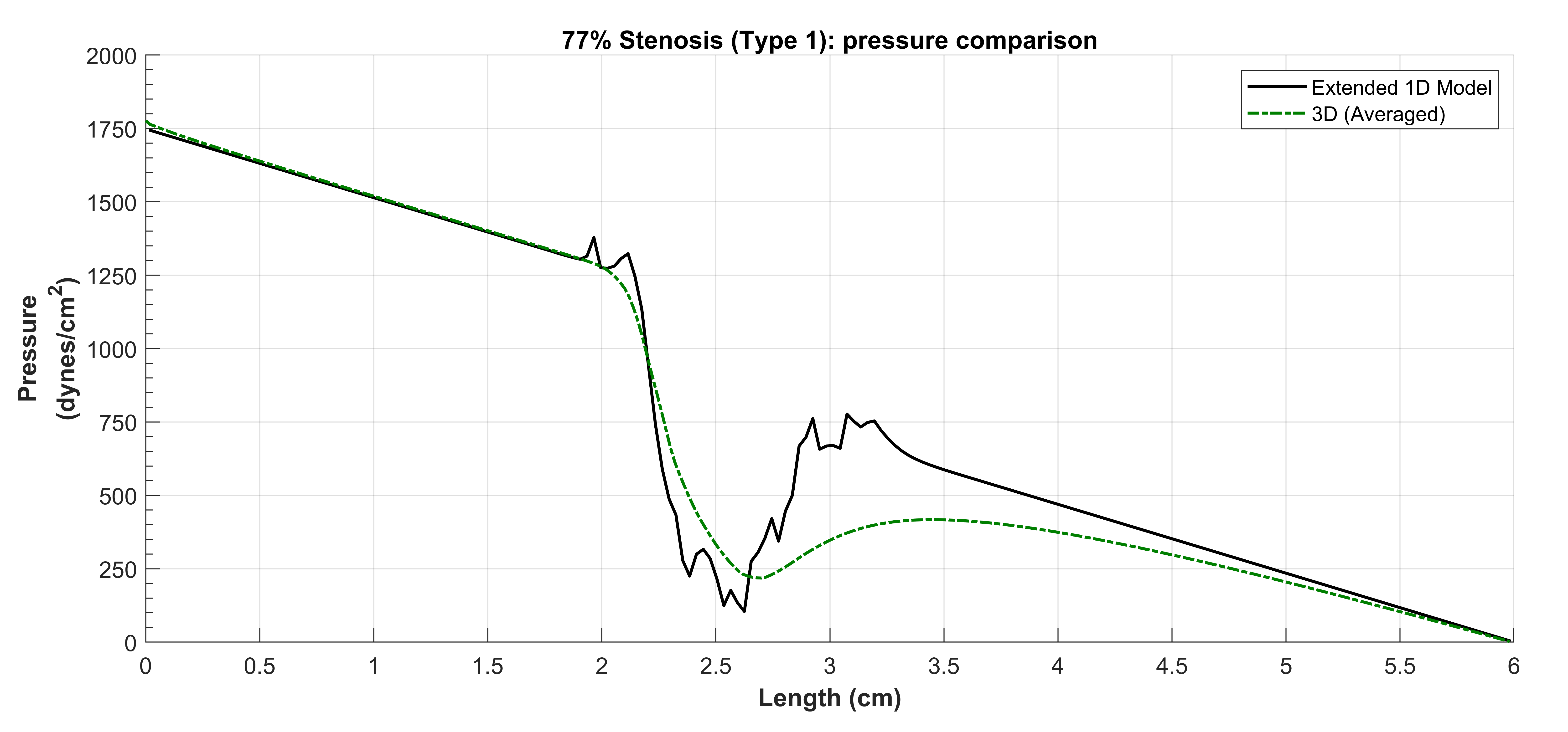}}
  \subfigure[40\% stenosis]{   \includegraphics[width=0.52\columnwidth]{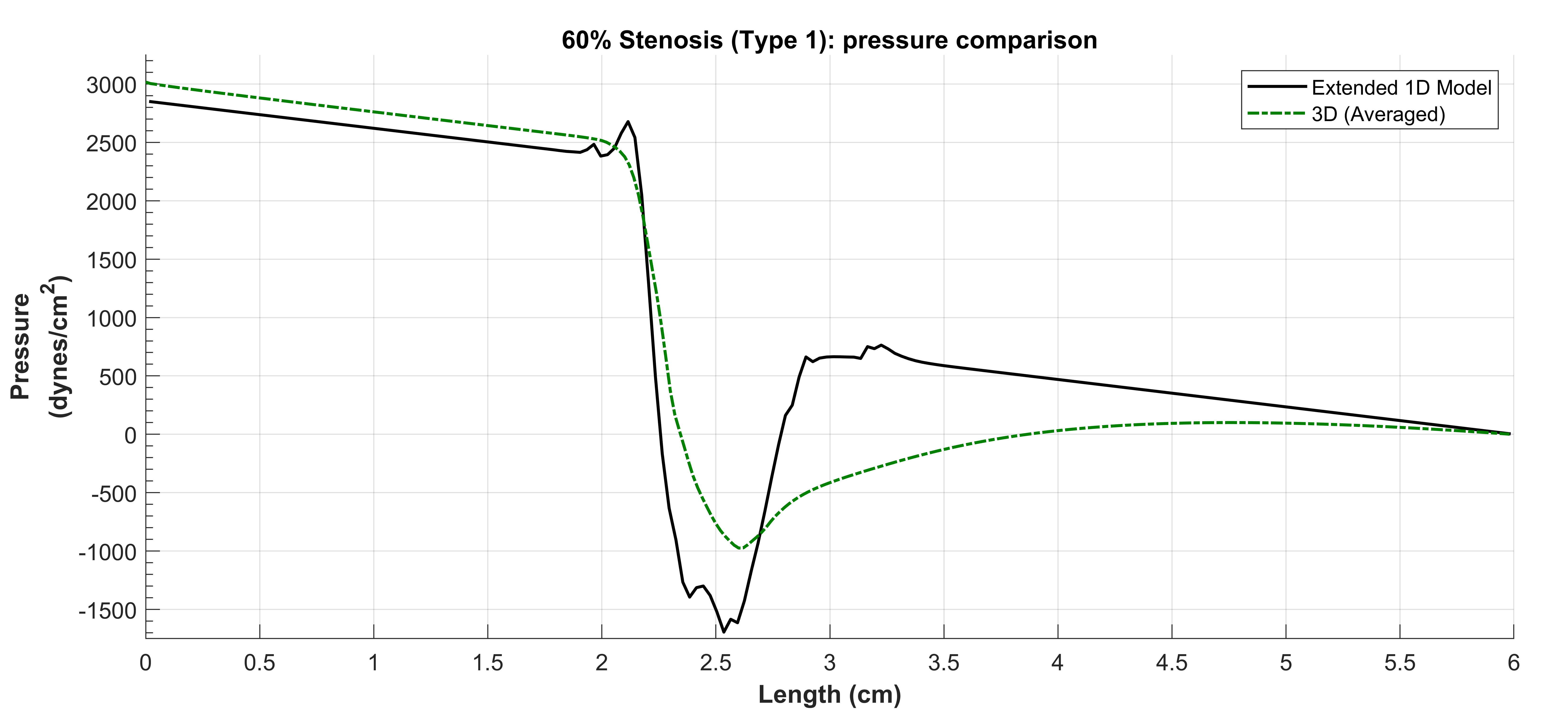}}
  \subfigure[50\% stenosis]{ \includegraphics[width=0.52\columnwidth]{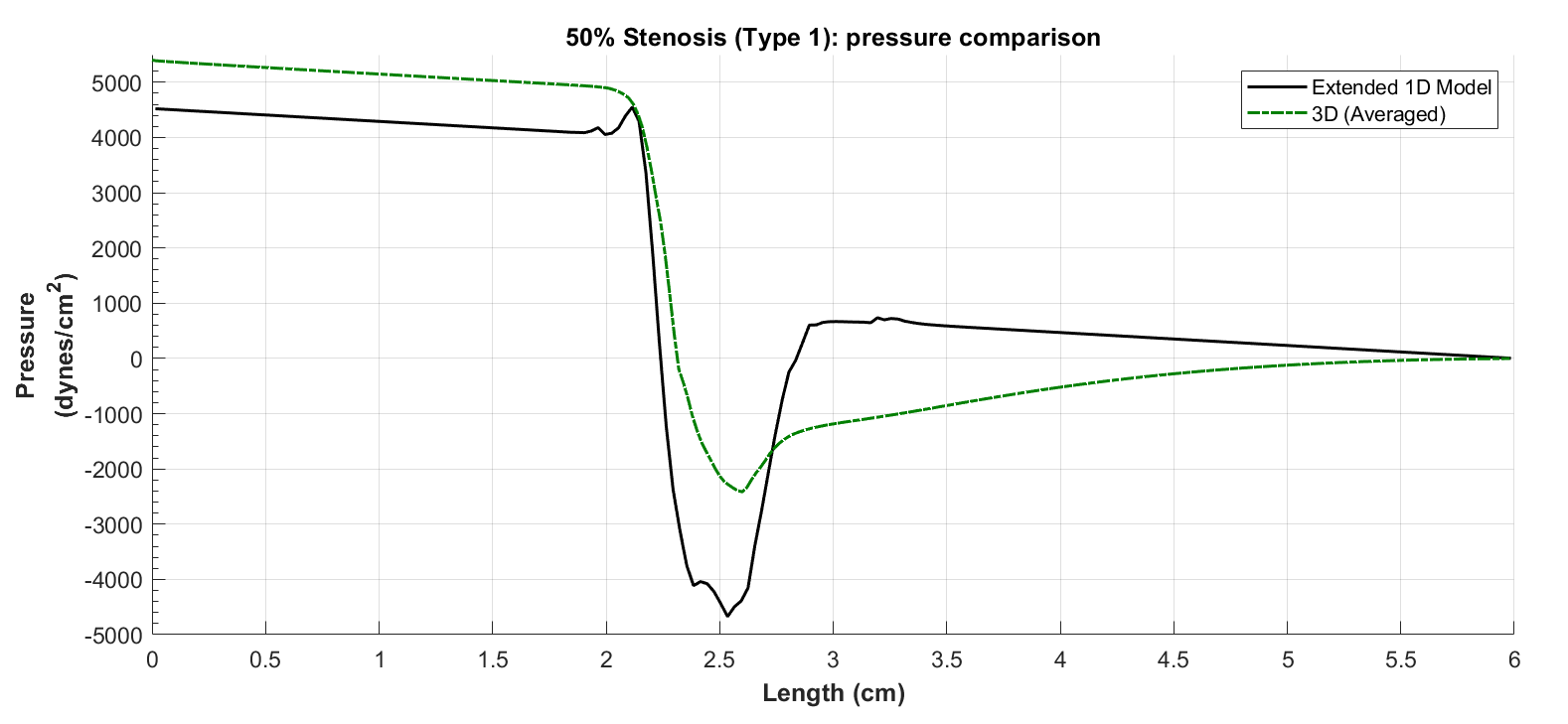}}
   \subfigure[pressure obtained from established $1D$ results]{ \includegraphics[width=0.52\columnwidth]{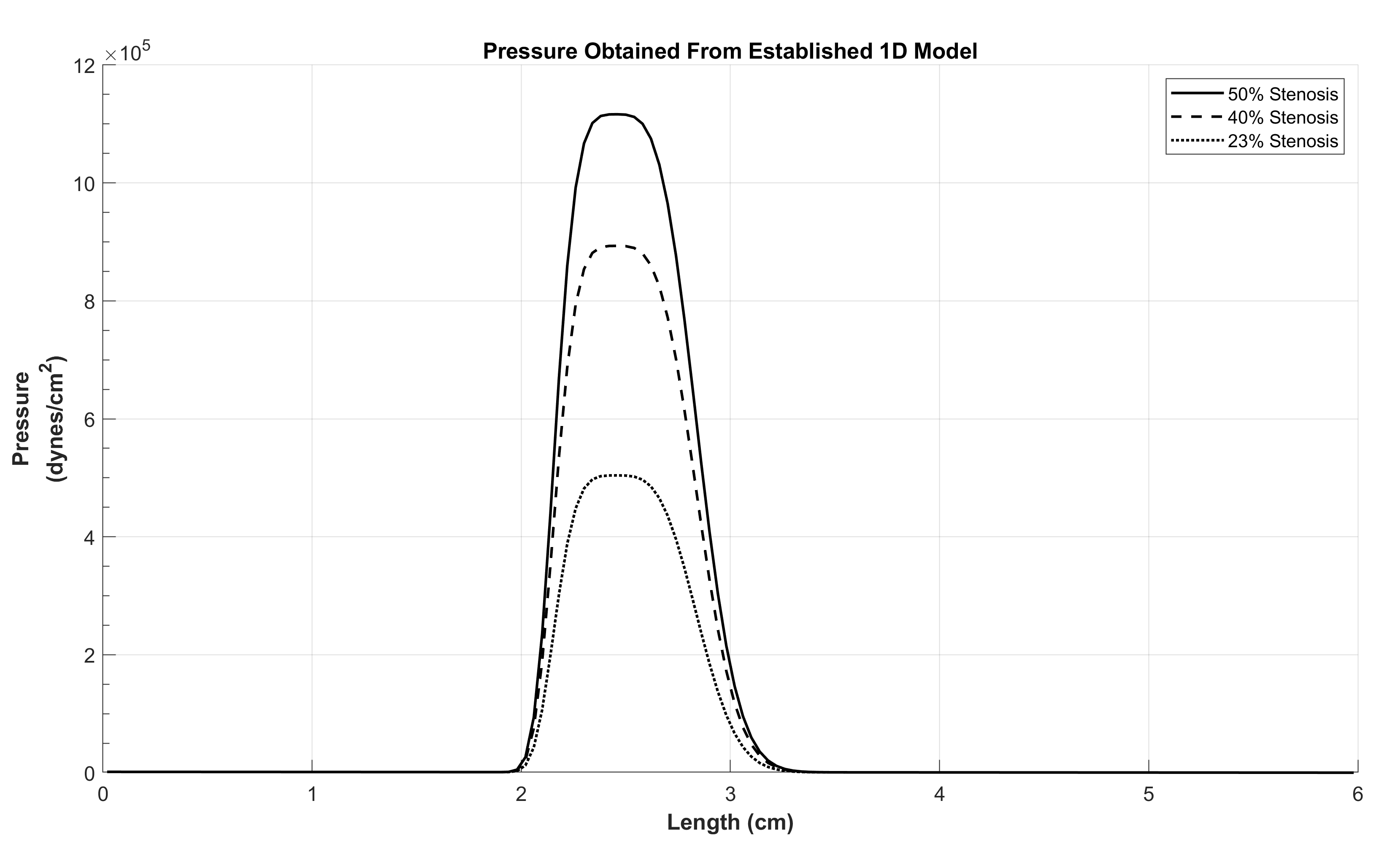}}
\caption{The pressure $p$ obtained from our extended $1D$ model (in green) is compared with the established $1D$ model (subfigure (d)) and the $3D$ full model (in black) results across all cases. Our extended $1D$ model again shows good agreement with the $3D$ results, while the established $1D$ model does not yield accurate solutions.}
\label{Fig.pressure}
\end{figure}

\if(0)
\begin{figure}[htbp]
	\center
	\subfigure[Type~1]{ \includegraphics[width=0.6\columnwidth]{./case1_ur.png}}
	\subfigure[Type~2]{   \includegraphics[width=0.6\columnwidth]{./case2_ur.png}}
	\subfigure[Type~3]{ \includegraphics[width=0.5\columnwidth]{./case3_ur.png}}
	\caption{Superimposed radial velocity $\eta_r$ of $1D$ and $3D$ results at proximal and distal for all three types}
	\label{Fig.ur}
\end{figure}
Next, we report the results of computed radial velocity in Figure \ref{Fig.ur}. Parabolic velocity profile as well as the velocity magnitude obtained from the $1D$ model agrees with the $3D$ results.
\fi

\begin{figure}[h]
	\centering
	\begin{overpic}[width=.7\textwidth,grid=false,tics=10]{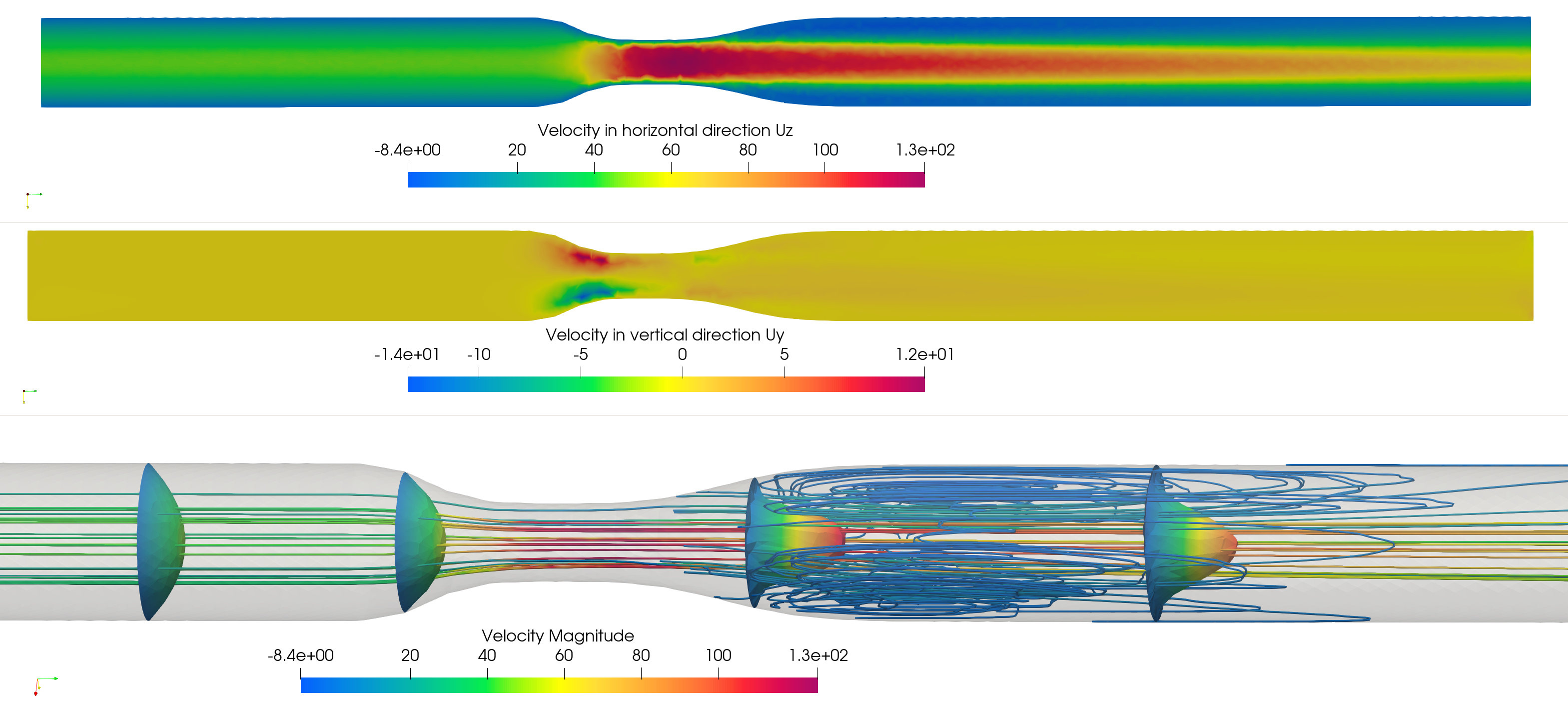}
	\end{overpic}
	\caption{Full $3D$ simulation results showing the velocity profile and streamlines}
	\label{3D_streamlines}
\end{figure}

Next, we present the pressure distribution along the artery obtained from our extended $1D$ model for all cases, overlaid with the results from the $3D$ simulation, as illustrated in Figure \ref{Fig.pressure}. Our extended $1D$ model effectively captures the overall trend of pressure changes in the longitudinal direction, except near the pre-stenotic region and post-stenotic area. In contrast, the established $1D$ model fails to relect the pressure trend. The mismatches between our extended $1D$ model results and $3D$ results in those areas attribute primarily the jet flow and recirculation effects, as illustrated in the Figure ~\ref{3D_streamlines}, which the extended $1D$ model fails to capture.

\if{0}
\subsection{Comparison between $1D$ reduced model and full $3D$ compliant vessel simulation with physiological data}
In this section, we would like to illustrate that our model is capable of handling physiological data
and providing meaningful results. We study the same stenosed coronary artery cases as in the previous section but using the physiological data from \cite{Marques2002}. The inlet and outlet pressure measurements as shown in Figure \ref{Fig.data} are prescribed as the boundary data in our $1D$ simulation.
In Figure \ref{wave_form}, we see that the fluid velocity wave forms at the stenotic and pre-stenotic region over one cardiac cycle for all three types. The numerical results agree well with the measurements reported in (\cite{Hozumi2000}). We notice that the velocity wave does not show significant difference between the three types of stenotic lesions, but it does increase significantly in the stenotic region comparing with the non-stenotic region.
\begin{figure}[h]
	\centering
	\begin{overpic}[width=0.7\textwidth,grid=false,tics=10]{data.png}
	\end{overpic}
	\caption{Inlet and outlet pressure data.}
	\label{Fig.data}
\end{figure}

\begin{figure}[htbp]
	\center
	\subfigure[Velocity at the stenosis]{ \includegraphics[width=0.8\columnwidth]{./VelocityAtStenosis.png}}
	\subfigure[Velocity at prestenosis]{ \includegraphics[width=0.8\columnwidth]{./VelocityAtPreStenosis.png}}
	\caption{Superimposed longitudinal velocity $u_z$ of $3D$ results for all cases.}
	\label{wave_form}
\end{figure}
\fi

\section{Conclusion}
We developed an extended $1D$ reduced model to analyze blood flow in arteries with stenosis. This model improves upon the well established $1D$ approach by incorporating the variable radius of the blood vessel, achieved through additional terms that account for the effects of the changing reference radius. The extended $1D$ model successfully reproduces the results and trends observed in the full $3D$ simulation, while the conventional $1D$ model fails to produce accurate outcomes. In addition, we introduce a post-processing technique that extracts radial velocity information from the $1D$ results, providing a more detailed understanding of the flow. Numerical simulations were performed for three different of stenotic lesions to demonstrate the accuracy and effectiveness of the proposed model.

In future work, we aim to extend the proposed model to study the blood flow dynamics in a full vascular network, encompassing more complex geometries and interactions within the circulatory system. By incorporating additional physiological factors, such as vessel elasticity and branching, we hope to further enhance the model's accuracy in simulating realistic blood flow conditions. Ultimately, we envision that this model could serve as a powerful tool for diagnostic purposes, aiding in the identification and assessment of vascular diseases such as stenosis, aneurysms, and other hemodynamic abnormalities. With further validation and refinement, the model has the potential to assist clinicians in making more informed decisions about treatment strategies.

\section{Acknowledgments}
This work has been supported in part by the following grants: NSF DMS-2247000 (Canic), NSF DMS-2247001 (Wang), Simon Foundation MP-TSM-00002663 (Wang).

\appendix
\section{Derivation of conservation of mass in cylindrical coordinates with axial symmetry and variable radius}
\begin{figure}[ht!]
\centering
\begin{overpic}[percent,grid=false,tics=20,scale=0.5]{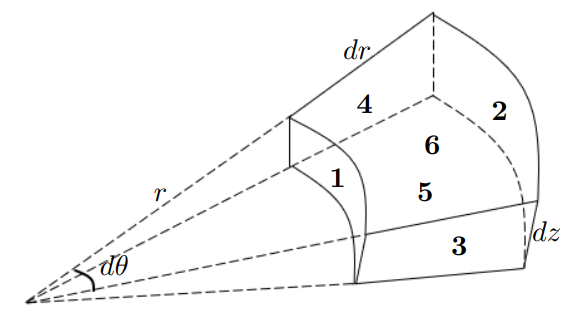}
\end{overpic}
	\caption{
	{Sketch of an infinitesimal volume inside a straight tube with variable radius}}
	\label{Sketch_cm}
\end{figure}

To formulate the equation for mass conservation in cylindrical coordinates, incorporating axial symmetry and variable radius, we commence by defining the alteration in mass. Considering within an infinitesimal volume $d V=r d \theta d r d z$ as illustrated in Figure \ref{Sketch_cm}, the change in mass is expressed as:
$$
\begin{aligned}
 \text { Change in mass }=\frac{d m}{d t}=\frac{d \rho V}{d t}.
&
\end{aligned}
$$
Next, we will examine the mass fluxes involving flow velocity components $u_r$, $u_\theta$, and $u_z$, as they enter and exit through each pair of faces.

1) Fluxes in the radial direction $\vec{r}$, representing the flow across faces 1 and 2 (Omitting higher order terms):
$$
\left\{~
\begin{aligned}
& \left(Flux_{\vec{r}}\right)_1=\underbrace{r d \theta d z}_{\text {Area of face } 1}\left(\rho u_r\right), \\
& \left(Flux_{\vec{r}}\right)_2=\underbrace{(r+d r) d \theta d z}_{\text {Area of face } 2}\left(\rho u_r+\frac{\partial \rho u_r}{\partial r} d r\right).
\end{aligned}
\right.
$$
$$
\begin{aligned}
\text { Total flux in } \vec{r} \text { direction }&=\left(Flux_{\vec{r}}\right)_2-\left(Flux_{\vec{r}}\right)_1 \\
& =r{d\theta dz}\left(\rho u_r\right)+ r d \theta d z\left(\frac{\partial \rho u_r}{\partial r} d r\right) + d r d \theta d z\left(\rho u_r\right)\\
&+d \theta d z \frac{\partial \rho u_r}{\partial r} d r^2 - r{d\theta dz}\left(\rho u_r\right) \\
& =r d \theta d z d r \frac{\partial \rho u_r}{\partial r}+d r d \theta d z\left(\rho u_r\right) \\
& =d V \frac{\partial \rho u_r}{\partial r}+\frac{d V}{r}\left(\rho u_r\right) =d V\left(\frac{\partial \rho u_r}{\partial r}+\frac{\rho u_r}{r}\right). 
\end{aligned}
$$

2) Fluxes in angular direction $\vec{\theta}$, representing the flow through faces 3 and 4:
$$
\left\{~
\begin{aligned}
& \left(Flux_{\vec{\theta}}\right)_3=\underbrace{d r d z}_{\text {Area of face } 3}\left(\rho u_\theta\right), \\
& \left(Flux_{\vec{\theta}}\right)_4=\underbrace{d r d z}_{\text {Area of face } 4}\left(\rho u_\theta+\frac{\partial \rho u_\theta}{\partial \theta} d \theta\right).
\end{aligned}
\right.
$$
$$
\begin{aligned}
\text { Total flux in } \vec{\theta} \text { direction }  =(Flux_{\vec{\theta}})_4-\left(Flux_{\vec{\theta}}\right)_3  =d r d z d \theta \frac{\partial \rho u_\theta}{\partial \theta}  =r d r d z d \theta \frac{1}{r} \frac{\partial \rho u_\theta}{\partial \theta}  =d V\left(\frac{1}{r} \frac{\partial \rho u_\theta}{\partial \theta}\right). \\
\end{aligned}
$$

3) Fluxes in axial direction $\vec{z}$, representing the flow across faces 5 and 6:
$$
\left\{~
\begin{aligned}
& \left(Flux_{\vec{z}}\right)_5=\underbrace{r d \theta d r}_{\text {Area of face } 5}\left(\rho u_z\right), \\
& \left(Flux_{\vec{z}}\right)_6=\underbrace{r d \theta d r}_{\text {Area of face } 6}\left(\rho u_z+\frac{\partial \rho u_z}{\partial z} d z\right). 
\end{aligned}
\right.
$$
$$
\text { Total flux in } \vec{z} \text { direction } =  \left(Flux_{\vec{z}}\right)_6-  \left(Flux_{\vec{z}}\right)_5 = r d \theta d r d z\left(\frac{\partial \rho u_z}{\partial z}\right)  =d V\left(\frac{\partial \rho u_z}{\partial z}\right).
$$
Thus, within the infinitesimal volume, we derive:
$$
d V \frac{d \rho}{d t}=d V\left(\frac{\partial \rho u_r}{\partial r}+\frac{\rho u_r}{r}+\frac{1}{r} \frac{\partial \rho u_\theta}{\partial \theta}+\frac{\partial \rho u_z}{\partial z}\right). 
$$
Assuming it holds everywhere, we obtain:
$$
\frac{d \rho}{d t}=\frac{\partial \rho u_r}{\partial r}+\frac{\rho u_r}{r}+\frac{1}{r} \frac{\partial \rho u_\theta}{\partial \theta}+\frac{\partial \rho u_z}{\partial z}.
$$
Given the assumption of incompressible fluid ($\rho$ is constant) and axial symmetry ($\frac{\partial u_\theta}{\partial \theta}=0$), we derive the equation for mass conservation as appears in Equ.~(\ref{NS}c). Noting that this equation is not in conservative form.

\section{Another option to approximate the integral term $\mathcal{I}$}
\label{section:AppendixB}
Here, we present an alternative approach to approximate the integral term $\mathcal{I}$. Referring back to Equ.~(\ref{approx}), we obtain:
$$
\int_0^{\tilde{R}} \tilde{r} \tilde{u}_r \tilde{u}_z d \tilde{r} =-2 \tilde{U} \frac{\partial \tilde{U}}{\partial \tilde{z}}\left(\frac{11}{105} \tilde{R}^3\right)-2 \tilde{U}^2 \left(\frac{\partial \ln \tilde{R}}{\partial \tilde{z}}
+\frac{\partial \ln R_0(\tilde{z})}{\partial \tilde{z}} \right)\left(\frac{2}{35} \tilde{R}^3\right).
$$
By choosing the following approximation $\frac{11}{105}\approx\frac{12}{105}=\frac{4}{35}$, one yields:
\begin{equation}
    \begin{aligned}
\mathcal{I}\approx&-\frac{\partial \ln R_0(\tilde{z})}{\partial \tilde{z}} 
\frac{1}{2 L^2} \frac{\partial R_0^2(\tilde{z})}{\partial \tilde{z}}
\left(\frac{4}{35}\tilde{R}^3\right)
\left[
\frac{\partial \tilde{U}^2}{\partial \tilde{z}}
+\tilde{U}^2\left(\frac{\partial ln R_0(\tilde{z})}{\partial \tilde{z}}
+\frac{\partial ln \tilde{R}}{\partial \tilde{z}}\right)
\right]\\
=&-\frac{4}{35L^2}\left(\frac{\partial R_0(\tilde{z})}{\partial \tilde{z}}\right)^2
\frac{\tilde{A}^2}{\tilde{R}}
\left[
\frac{\partial \tilde{U}^2}{\partial \tilde{z}}
+\tilde{U}^2\left(\frac{\partial ln R_0(\tilde{z})}{\partial \tilde{z}}
+\frac{\partial ln \tilde{R}}{\partial \tilde{z}}\right)
\right]\\
=&-\frac{4}{35L^2}\left(\frac{\partial R_0(\tilde{z})}{\partial \tilde{z}}\right)^2
\frac{1}{\tilde{R}}
\left[
\frac{\partial \tilde{Q}^2}{\partial \tilde{z}}
-\tilde{U}^2\frac{\partial \tilde{A}^2}{\partial \tilde{z}}
+\tilde{Q}^2\left(\frac{\partial ln R_0(\tilde{z})}{\partial \tilde{z}}
+\frac{\partial ln \tilde{R}}{\partial \tilde{z}}\right)
\right].
    \end{aligned}
\end{equation}
After converting the above equation into its dimensional form, we have the expression for term $\mathcal{I}$ as follows:
$$
\begin{aligned}
\mathcal{I}=&-\frac{4}{35}\left(\frac{\partial R_0}{\partial z}\right)^2
    \frac{R_0}{R}
\left[\left(
\frac{L}{R_0^4 U_z^2}\frac{\partial Q^2}{\partial z}
-\frac{4 L Q^2}{R_0^4 U_z^2}\frac{\partial ln R_0}{\partial z}\right)
-\left(\frac{U^2 L}{U_z^2 R_0^4}\frac{\partial A^2}{\partial z}
-\frac{4 U^2 L A^2}{U_z^2 R_0^4}\frac{\partial ln R_0}{\partial z}
\right)\right.\\
&\left.+\frac{Q^2 L}{R_0^4 U_z^2}\frac{\partial ln R}{\partial z}
\right]\\
=&-\frac{4}{35}\left(\frac{\partial R_0}{\partial z}\right)^2
    \frac{R_0}{R}
\left[
\frac{L}{R_0^4 U_z^2}\frac{\partial Q^2}{\partial z}
-\frac{U^2 L}{U_z^2 R_0^4}\frac{\partial A^2}{\partial z}
+\frac{Q^2 L}{R_0^4 U_z^2}\frac{\partial ln R}{\partial z}
\right]\\
=&-\frac{4}{35}\left(\frac{\partial R_0}{\partial z}\right)^2
    \frac{1}{R R_0}
 \left(\frac{L}{U_z^2 R_0^2} \right)
\left[
\frac{\partial Q^2}{\partial z}
-U^2 \frac{\partial A^2}{\partial z}
+Q^2\frac{\partial ln R}{\partial z}
\right]\\
=&-\frac{4}{35}\left(\frac{\partial R_0}{\partial z}\right)^2
    \frac{1}{RR_0}\left(\frac{L}{U_z^2 R_0^2} \right)\left(\frac{\partial Q^2}{\partial z}-\frac{Q^2}{2A^2}\frac{\partial A^2}{\partial z}\right)\\
=&-\frac{4}{35}\left(\frac{\partial R_0}{\partial z}\right)^2
    \frac{1}{R_0}\left(\frac{L}{U_z^2 R_0^2} \right)\left(\frac{1}{\sqrt{A}}\frac{\partial Q^2}{\partial z}-\frac{Q^2}{2A^{3/2}}\frac{\partial A^2}{\partial z}\right)\\
=&-\frac{4}{35}\left(\frac{\partial R_0}{\partial z}\right)^2
    \frac{1}{R_0}\left(\frac{L}{U_z^2 R_0^2} \right)\frac{\partial}{\partial z}\left(\frac{Q^2}{\sqrt{A}}\right).
\end{aligned}
$$
Therefore, another version of the rectified (A, Q) system is given as follows:
\begin{equation}
\left\{	
\begin{aligned}
		&\frac{\partial A}{\partial t}
		+ \frac{\partial Q}{\partial z} = 0,\\
		&\frac{\partial Q}{\partial t}
		+ \frac{\partial}{\partial z}\left(\alpha \frac{Q^{2}}{A }\right)
		+\frac{A}{\rho_{f}} \frac{\partial p}{\partial z}
\textcolor{black}{-\frac{4}{35}\left(\frac{\partial R_0}{\partial z}\right)^2\frac{1}{R_0}\frac{\partial}{\partial z}\left(\frac{Q^2}{\sqrt{A}}\right)}
		=\frac{2}{Re} R\frac{\d u_z(R, z, t)}{\d r}.
	\end{aligned}
 \right.
\end{equation}

\bibliographystyle{elsarticle-num} 

\bibliography{reference}

\end{document}